\documentclass[11pt]{article}
\usepackage{geometry,amsmath,amssymb,physics,booktabs,graphicx,caption,authblk}
\usepackage{biblatex}
\usepackage{xcolor}
\usepackage[colorlinks, linkcolor=blue, citecolor=black!50!green]{hyperref}
\usepackage[thmmarks,amsmath,thref]{ntheorem}
\theorembodyfont{\normalfont}
\newtheorem{theorem}{Theorem}
\newtheorem{lemma}{Lemma}
\newtheorem{definition}{Definition}

\newtheorem{remark}{Remark}
\newcommand{\setzero}{\setlength{\itemsep}{0pt}}
\renewcommand{\d}{\mathrm{d}}

\newcommand{\T}{\mathrm{T}}
\newcommand{\Ge}{\geqslant}
\newcommand{\Le}{\leqslant}
\newcommand{\Succ}{\succcurlyeq}
\newcommand{\Prec}{\preccurlyeq}
\title{Statistical Error of Numerical Integrators for Underdamped Langevin Dynamics with Deterministic And Stochastic Gradients}
\author[1]{Xuda Ye%
\thanks{abneryepku@pku.edu.cn}}
\author[2]{Zhennan Zhou%
\thanks{zhouzhennan@westlake.edu.cn}}
\affil[1]{Beijing International Center for Mathematical Research, Peking University, Beijing, 100871, China.}
\affil[2]{Institute for Theoretical Sciences, Westlake University, Hangzhou, 310030, China.}

\begin{document}
%% ---- body text ---- %%
\maketitle

\begin{abstract}
We propose a novel discrete Poisson equation approach to estimate the statistical error of a broad class of numerical integrators for the underdamped Langevin dynamics.
The statistical error refers to the mean square error of the estimator to the exact ensemble average with a finite number of iterations. With the proposed error analysis framework, we show that when the potential function $U(x)$ is strongly convex in $\mathbb R^d$ and the numerical integrator has strong order $p$, the statistical error is $\mathcal O(h^{2p}\hspace{-1pt}+\hspace{-1pt}\frac1{Nh})$, where $h$ is the time step and $N$ is the number of iterations. Besides, this approach can be adopted to analyze integrators with stochastic gradients, and quantitative estimates can be derived as well.
Our approach only requires the geometric ergodicity of the continuous-time underdamped Langevin dynamics, and relaxes the constraint on the time step.
\end{abstract}

\noindent
\textbf{Keywords} underdamped Langevin dynamics, geometric ergodicity, discrete Poisson equation, statistical error\\[6pt]
\textbf{MSC2020} 60H35, 37M05

\section{Introduction}
\label{section: 1}

The computation of the ensemble average from a target Boltzmann distribution stands as a cornerstone challenge in both molecular dynamics \cite{md_1,BAOAB_1,Poisson_4} and data science \cite{SGLD_1,SGLD_0}. We examine the Hamiltonian system
\begin{equation*}
	H(x,v) = \frac{|v|^2}2 + U(x),
\end{equation*}
where $x$ and $v$ represent the particle's position and velocity in $\mathbb R^d$, respectively, and $U(x)$ denotes the potential function in $\mathbb{R}^d$.
Our primary focus is the error analysis of numerical methods employed for sampling from the target Boltzmann distribution
\begin{equation*}
	\pi(x,v) = \frac1Z e^{-H(x,v)},~~~~
	Z = \int_{\mathbb R^{2d}}
	e^{-H(x,v)}\d x\d v.
\end{equation*}

In the realm of molecular dynamics, an arsenal of numerical methods has been proposed with the aim of effectively sampling from $\pi(x,v)$. These methods encompass Langevin dynamics \cite{md_1,BAOAB_1,md_2}, Hamiltonian Monte Carlo (HMC) \cite{HMC_1,HMC_2}, Andersen dynamics \cite{A_1,A_2}, and Nosé–Hoover dynamics \cite{N_1}. In this paper, we focus on the time discretization of the underdamped Langevin dynamics
\begin{equation}
\left\{
\begin{aligned}
\dot x_t & = v_t, \\
\dot v_t & = -\nabla U(x_t) - \gamma v_t + \sqrt{2\gamma}\dot B_t,
\end{aligned}
\right.
\label{underdamped}
\end{equation}
where $x_t$ and $v_t$ denote the position and velocity variables in $\mathbb{R}^d$, respectively. Here, $\gamma>0$ represents the damping rate, and $(B_t)_{t\Ge 0}$ signifies the Brownian motion in $\mathbb{R}^d$. By fixing the test function $f(x,v)$ in $\mathbb R^{2d}$, the ergodic theory dictates that the ensemble average $\pi(f)$ can be obtained through the limit of the time average
\begin{equation*}
	\pi(f) := \int_{\mathbb R^{2d}}
	f(x,v) \pi(x,v)\d x\d v = 
	\lim_{T\rightarrow\infty}
	\frac1T \int_0^T f(x_t,v_t)\d t. 
\end{equation*}

The time discretization of the underdamped Langevin dynamics \eqref{underdamped} lends itself to implementation via a numerical integrator, offering a stochastic approximation to the ensemble average $\pi(f)$ attainable in practice. Suppose a numerical integrator with time step $h$ yields the numerical solution $(X_n,V_n)_{n\Ge 0}$; in that case, the error in computing $\pi(f)$ can be formulated as
\begin{equation}
	e(N,h) = \frac1N \sum_{n=0}^{N-1} f(X_n,V_n) - \pi(f),
	\label{time average error}
\end{equation}
where $N\in\mathbb{N}$ denotes the number of iterations.
As $N$ tends to infinity and $h$ decreases to 0, the error $e(N,h)$ is expected to converge to 0, as the consistency of a numerical integrator requires. Therefore,
the accuracy of the numerical integrator can be characterized by:
\begin{itemize}
\setzero
\item the \emph{bias} $\mathbb{E}[e(N,h)]$, delineating the deviation of the numerical invariant distribution;
\item the \emph{statistical error} $\mathbb{E}[e^2(N,h)]$, emanating from the bias and the random fluctuation of the numerical solutions $(X_n,V_n)_{n\Ge0}$, also referred to as the \emph{mean square error}.
\end{itemize}
While the bias indicates the average error of an ensemble of stochastic numerical solutions, the statistical error encapsulates the accuracy of the a single trajectory in a more holistic manner.
The objective of this paper is to estimate the statistical error of various numerical integrators for the underdamped Langevin dynamics \eqref{underdamped}.

We study a large class of prevailing numerical integrators such as the Euler--Maruyama \eqref{EM} integrator and the UBU \eqref{UBU} integrator \cite{UBU_1,UBU_2,UBU_3}, also their stochastic gradient versions. It is noteworthy that \eqref{EM} possess a strong order of 1, while \eqref{UBU} boasts a strong order of 2. Below, we informally outline the main findings of this paper. The full statement of the theorems can be found in Sections~\ref{section: 5} \& \ref{section: 6}.\\[6pt]
\textbf{Theorem~\ref{theorem: numerical error}.}
Suppose the numerical method has strong order $p$. If $U(x)$ is \emph{strongly convex outside a ball}, then there exists a constant $\gamma_0>0$ such that when the damping rate $\gamma \Ge \gamma_0$, 
\begin{equation*}
	\mathbb E[e^2(N,h)] = \mathcal O\bigg(
	h^{2p-1} + \frac1{Nh}
	\bigg).
\end{equation*}
\textbf{Theorem~\ref{theorem: numerical error convex}.}
Suppose the numerical method has strong order $p$. If $U(x)$ is \emph{strongly convex in $\mathbb R^d$}, then there exists a constant $\gamma_0>0$ such that when the damping rate $\gamma \Ge \gamma_0$,
\begin{equation*}
	\mathbb E[e^2(N,h)] = \mathcal O\bigg(
	h^{2p} + \frac1{Nh}
	\bigg).
\end{equation*}
When $U(x)$ is strongly convex in $\mathbb R^d$, we study the statistical error of the stochastic gradient numerical integrators, including the stochastic gradient Euler--Maruyama \eqref{SG-EM} and the stochastic gradient UBU \eqref{SG-UBU}.\\[6pt]
\textbf{Theorem~\ref{theorem: SG}.}
Suppose the stochastic gradient numerical integrator has strong order $p$.
If $U(x)$ is \emph{strongly convex in $\mathbb R^d$}, then there exists a constant $\gamma_0>0$ such that when the damping rate $\gamma \Ge \gamma_0$,
\begin{equation*}
\mathbb E[e^2(N,h)] = \mathcal O\bigg(
h^{\min\{2p,2\}} + \frac1{Nh}
\bigg).
\end{equation*}

Before we illustrate the main strategy to prove these results, we review the prevailing approaches to estimate the statistical error. In general, the bias of the numerical method can be estimated via the geometric ergodicity of the numerical method \cite{Lyapunov,error_1,error_2}, namely, the convergence rate of the distribution law of the numerical solution.
In comparison, the statistical error usually encompasses more information than the geometric ergodicity, and the analysis based on the functional inequalities \cite{left_2,LS_1,LS_2,LS_3} can be inappropriate for estimating the statistical error.
Currently, there are two major approaches to analyze the statistical error for the underdamped Langevin dynamics \eqref{underdamped}. 

The first one is to prove the global contractivity of the numerical integrator. 
The utilization of the probabilistic coupling approach facilitates the establishment of global contractivity for a broad spectrum of numerical integrators, including both HMC \cite{coupling_2,coupling_3,coupling_6,coupling_7} and Langevin discretization \cite{coupling_1,coupling_4,coupling_5}. Subsequently, leveraging these contractivity results allows for a direct and explicit estimation of the statistical error.
However, the application of these results necessitates stringent conditions on the parameter: the time step $h$ should be sufficiently small to ensure the global contractivity in a single step. When the numerical integrator is driven by stochastic gradients, it's important to note that these contractivity results not always guarantee the optimal order in the time step \cite{coupling_2,coupling_6}.

The second one is to utilize the Poisson equation.
The Poisson equation approach has gained prominence in the analysis of overdamped Langevin dynamics since the seminal work \cite{Poisson_1}. Subsequently, this approach has been utilized to scrutinize the statistical error of the Euler--Maruyama integrator \cite{Poisson_2,Poisson_3} and the stochastic gradient Langevin dynamics (SGLD) \cite{SGLD_1,SGLD_2,SGLD_3,SGLD_4,SGLD_5}. Notably, the Poisson equation approach shares a close connection with Stein's method \cite{Poisson_1,Poisson_3}, and serves as a viable strategy to numerically augment the sampling efficiency of underdamped Langevin dynamics \cite{Pavliotis_1,Pavliotis_2,Pavliotis_3}.
It's noteworthy that the efficacy of the Poisson equation approach primarily hinges upon the ability to pointwise estimate the solution to the Poisson equation \cite{Poisson_3}. In the overdamped scenario, such estimates can be achieved via $L^p$-estimates, owing to the fact that the generator $\mathcal{L}$ is a uniformly elliptic operator. However, in the underdamped case, such estimates can pose significant challenges due to the hypoellipticity of the generator $\mathcal{L}$, thereby rendering the estimation of statistical error via the Poisson equation approach inherently difficult.

In this paper, we estimate the statistical error of numerical integrators for the underdamped Langevin dynamics \eqref{underdamped} utilizing the discrete Poisson solution $\phi_h(x,v)$. This solution $\phi_h(x,v)$ can be perceived as a time-discretized version of the original Poisson solution $\phi(x,v)$ (see Section~\ref{section: 5}). The utilization of $\phi_h(x,v)$ enables the explicit representation of the time average in \eqref{rountine 2}, and thus circumvents the need for high-order Taylor expansions on $\phi(x,v)$. Consequently, the statistical error of a wide class of numerical integrators can be examined using the same discrete Poisson solution $\phi_h(x,v)$.
Moreover, the constraint on the time step $h$ is solely dependent on the uniform-in-time moment estimate in Lemmas \ref{lemma: moment}, which is typically less stringent than the requirement in the proof of global contractivity.
When the potential function $U(x)$ is strongly convex in $\mathbb{R}^d$, we introduce a novel estimation approach for the second derivative $\nabla^2 \phi_h(x,v)$ based on the tangent process \cite{Pavliotis_2}. This estimate of $\nabla^2\phi_h(x,v)$ facilitates the examination of the statistical error of stochastic gradient numerical integrators, via the explicit representation of the time average in \eqref{e SG}. 

We present a brief comparison of our results and those of related works. In the case of deterministic gradients, we provide the estimate of the statistical error for a broad class of numerical integrators in Theorems~\ref{theorem: numerical error} and \ref{theorem: numerical error convex}. When $U(x)$ is strongly convex in $\mathbb R^d$, the error order in the time step coincides with the strong order of the numerical integrator. Compared to those works based on the global contractivity of the numerical integrators \cite{UBU_3,coupling_1,coupling_2,coupling_4}, our analysis only requires the geometric ergodicity of the continuous-time dynamics \eqref{underdamped}, and the numerical stability is simply guaranteed by the uniform-in-time moments of the numerical solution.
In this way, we relax the constraint of the upper bound on the time step $h$. However, there are some limitations in our results for numerical integrators with deterministic gradients. First, we are unable to obtain the optimal statistical error $\mathcal O(h^{2p}\hspace{-1pt}+\hspace{-1pt}\frac1{Nh})$ when $U(x)$ is not globally convex. As a comparison, it is proved in \cite{coupling_3,coupling_5} that the numerical invariant distribution of the generalized Hamiltonian Monte Carlo (gHMC) has an $\mathcal O(h^2)$ deviation from the true invariant distribution $\pi$, and thus the statistical error is $\mathcal O(h^4\hspace{-1pt}+\hspace{-1pt}\frac1{Nh})$.
Also, the dependence on the parameters in Theorems~\ref{theorem: numerical error} and \ref{theorem: numerical error convex} are not clear, especially on the spatial dimension $d$.
Finally, since the error order in the time step $h$ is limited by the strong order, our results are suboptimal for the integrators whose weak order is larger than the strong order, such as the popular BAOAB integrator \cite{BAOAB_1,BAOAB_2,BAOAB_3}.

For numerical integrators with stochastic gradients, the proposed approach leads to improved results which are not justified before. In particular, Theorem~\ref{theorem: SG} showcases the optimal statistical error $\mathcal O(h^2\hspace{-1pt}+\hspace{-1pt}\frac1{Nh})$ when $U(x)$ is strongly convex in $\mathbb R^d$, for both \eqref{SG-EM} and \eqref{SG-UBU}. As a comparison, although the global contractivity of the stochastic gradient generalized Hamiltonian Monte Carlo (SGgHMC) has been well established in \cite{coupling_2,coupling_6}, these contractivity results do not provide an optimal error bound in terms of the stochastic gradients.
For instance, it is shown in \cite{coupling_2,coupling_6} that the deviation of the numerical invariant distribution of SGgHMC from the true invariant distribution $\pi$ is $\mathcal O(\sqrt{h})$, resulting the statistical error to be only $\mathcal O(h\hspace{-1pt}+\hspace{-1pt}\frac1{Nh})$.
Theorem~\ref{theorem: SG} also shows that there is no gain in the error order by using a second-order scheme for numerical integration, which aligns with the discourse presented in Proposition 18 of \cite{coupling_2}.

This paper is structured as follows.
In Section~\ref{section: 2}, we provide an introduction to the notations, assumptions, and the main idea of the proof. Specifically, we elucidate how the discrete Poisson solution $\phi_h(x,v)$ is utilized to represent the time average error.
Section~\ref{section: 3} introduces a specific definition of the strong order of the numerical integrator and establishes the uniform-in-time estimate.
In Section~\ref{section: 4}, we establish the estimates of the Kolmogorov equation, which serves as a prerequisite for estimating the discrete Poisson solution.
Section~\ref{section: 5} presents the estimation of numerical integrators with deterministic gradients.
Section~\ref{section: 6} delves into the estimation of numerical integrators with stochastic gradients.
In Section~\ref{section: 7}, we implement one numerical experiment to validate the convergence order of the statistical error for various numerical integrators.
\section{Setup and main ideas}
\label{section: 2}
We provide a list of notations, assumptions, and layout the main idea underlying the estimation of the statistical error.
\subsection{Notations}
The stochastic processes listed below are all in $\mathbb R^d\times\mathbb R^d$, namely, the position variable $x$ and velocity variable $v$ are in $\mathbb R^d$.\\[6pt]
\textbf{Underdamped Langevin dynamics.}
\begin{itemize}
\setzero
\item $z_t = (x_t,v_t)$ is the exact solution to the underdamped Langevin dynamics \eqref{underdamped} with the initial state $z_0 = (x_0,v_0)$.
\end{itemize}
\textbf{Numerical integrators with deterministic gradients.}
\begin{itemize}
\setzero
\item $Z_n = (X_n,V_n)$ is the numerical solution to \eqref{underdamped} with the initial condition $Z_0 = (X_0,V_0)$ and the time step $h$.
\item $Z_n(t) = (X_n(t), V_n(t))$ is the exact solution to \eqref{underdamped} with the initial value $Z_n = (X_n,V_n)$.
\end{itemize}
In this way, $Z_{n+1} - Z_n(h)$ represents the local error of the numerical integrator at each time step.\\[6pt]
\textbf{Stochastic gradient numerical integrators.}\\[6pt]
Let $b(x,\omega):\mathbb{R}^d\times\Omega\rightarrow\mathbb{R}^d$ be the stochastic gradient, where $\omega$ represents a random variable on the set $\Omega$. 
The random variables $(\omega_n)_{n\Ge0}$ are sampled independently to generate a numerical solution driven by stochastic gradients.
\begin{itemize}
\setzero
\item $\tilde Z_n = (\tilde X_n,\tilde V_n)$ is the numerical solution to the Langevin dynamics driven by stochastic gradients with the initial state $\tilde Z_0 = (\tilde X_0,\tilde V_0)$ and the time step $h$.
\item $\tilde Z_{n,1} = (\tilde X_{n,1},\tilde V_{n,1})$ is the one-step update of the integrator driven by the full gradient $\nabla U(x)$ with the initial state $\tilde Z_n = (\tilde X_n,\tilde V_n)$ and the time step $h$.
\item $\tilde Z_n(t) = (\tilde X_n(t),\tilde V_n(t))$ is the exact solution to \eqref{underdamped} with the initial state $\tilde Z_n = (\tilde X_n,\tilde V_n)$.
\item $\tilde Z_n\{t\} = (\tilde X_n\{t\},\tilde V_n\{t\})$ is the exact solution to the Langevin dynamics driven by the stochastic gradient $b(x,\omega_n)$ with the initial state $\tilde Z_n = (\tilde X_n,\tilde V_n)$.
\end{itemize}
Consequently, $\tilde Z_{n+1} - \tilde Z_{n,1}$ encapsulates the error due to stochastic gradients, and $\tilde Z_{n,1} - Z_n(h)$ represents the local error of the integrator with full gradients at each time step.

Given a fixed initial distribution, the exact solution $(x_t,v_t)_{t\Ge0}$ and the numerical solution $(X_n,V_n)_{n\Ge0}$ are determined by the Brownian motion $(B_t)_{t\Ge0}$, while the numerical solution with stochastic gradients $(\tilde X_n,\tilde V_n)_{n\Ge0}$ is determined by the Brownian motion $(B_t)_{t\Ge0}$ and the random variables $(\omega_n)_{n\Ge0}$ at each time step.\\[6pt]
\textbf{Universal constants.}
The universal constants $C$ and $\lambda$ are consistently employed throughout this paper. Their values do not depend on the time step $h$ and the number of iterations $N$, and may vary in different contexts. However, it's noteworthy that $C$ and $\lambda$ can still depend on other parameters such as the damping rate $\gamma$ and the dimension $d$.\\[6pt]
\textbf{Vector and matrices.}
For a vector, matrix or tensor $v$, the symbol $|v|$ means the Frobenius norm, namely, the Euclidean norm of all entries of $v$.
The symbols $O_d$ and $I_d$ represent the zero and identity matrices in $\mathbb R^{d\times d}$. The symbol $\Succ$ and $\Prec$ signify the Loewner order in $\mathbb R^{d\times d}$, which means $A\Succ B$ if and only if $A-B$ is positive semidefinite.
For a two-variable funtion $g(x,v)$, the symbols $\nabla_x g,\nabla_v g\in\mathbb R^d$ mean the gradients in $x$ and $v$ variables, respectively, and $\nabla g = (\nabla_x g, \nabla_v g) \in \mathbb R^{2d}$ is the full gradient.
\subsection{Assumptions}
The assumptions regarding the potential function $U(x)$ and the test function $f(x,v)$ are enumerated as follows. \\[6pt]
\textbf{A1~~}The potential function $U(x)\in C^3(\mathbb R^d)$, and for some constant $C_1>0$,
\begin{equation*}
|U(x)| \Le C_1(|x|+1)^2,~~~
|\nabla U(x)| \Le C_1(|x|+1),~~~
|\nabla^2 U(x)| \Le C_1,~~~
|\nabla^3 U(x)| \Le C_1.
\end{equation*}
\textbf{A2~~}$U(x)$ is \emph{strongly convex outside a ball}, namely for some constants $m,R>0$,
\begin{equation*}
\nabla^2 U(x) \Succ m I_d,~~~~\mathrm{for~any~}|x|\Ge R.
\end{equation*}
\textbf{A2'~~}$U(x)$ in \emph{strongly convex in $\mathbb R^d$}, namely for some constant $m>0$,
\begin{equation*}
\nabla^2 U(x) \Succ m I_d,~~~~\mathrm{for~any~}x\in\mathbb R^d.
\end{equation*}
\textbf{A3~~}The test function $f(x,v) \in C^3(\mathbb R^{2d})$, and for some constant $C_2>0$,
\begin{equation*}
|\nabla f(x,v)| \Le C_2,~~~~
|\nabla^2 f(x,v)| \Le C_2,~~~~
|\nabla^3 f(x,v)| \Le C_2.
\end{equation*}
% assumption on U(x)
The boundedness of $\nabla^2 U(x)$ and $\nabla^3 U(x)$ in \textbf{A1} is required for estimating the local error of the numerical integrators, as delineated in Lemma~\ref{lemma: strong order}.
% assumption on f(x)
The boundedness of $\nabla^2 f(x)$ and $\nabla^3 f(x)$ in \textbf{A3} is needed for estimating the Kolmogorov equation, as detailed in Lemma~\ref{lemma: qp diff decay}.
% explanations on assumptions
\textbf{A2} and \textbf{A2'} correspond to convexity outside the ball and global convexity in $\mathbb{R}^d$, respectively. Furthermore, \textbf{A2} entails the following drift condition:
\begin{equation}
x\cdot \nabla U(x) \Ge m |x|^2 - C_0,~~~~
\mathrm{for~any~}x\in\mathbb R^d,
\label{U diffu}
\end{equation}
where $C_0>0$ is a universal constant depending on $m,d,C_1$.

Lastly, we introduce the assumption regarding the stochastic gradient $b(x,\omega)$.\\[6pt]
\textbf{A4~~}For some constants $m,C>0$, $b(x,\omega)$ has linear growth,
\begin{equation*}
	|b(x,\omega)| \Le C_1(1+|x|),~~~~
	|\nabla b(x,\omega)| \Le C_1,
\end{equation*}
and $b(x,\omega)$ satisfy the drift condition
\begin{equation*}
	x\cdot \nabla b(x,\omega) \Ge m|x|^2 - C_0,~~~~
	\forall \omega\in\Omega.
\end{equation*}
Furthermore, $b(x,\omega)$ is an unbiased approximation to the gradient $\nabla U(x)$, namely
\begin{equation*}
\mathbb E^\omega [b(x,\omega)] = \nabla U(x).
\end{equation*}
The drift condition of $b(x,\omega)$ is employed in deriving the uniform-in-time moment bound of stochastic gradient numerical integrators, as expounded in Lemma~\ref{lemma: moment SG}.
\subsection{Main idea of the statistical error estimate}
Here's a concise procedure outlining how to estimate the statistical error $\mathbb E[e^2(N,h)]$ for various numerical integrators.
\begin{enumerate}
\item 
Given the test function $f(x,v)$, construct the Kolmogorov solution as
\begin{equation*}
	u(x,v,t) = (e^{t\mathcal L}f)(x,v) - \pi(f) = 
	\mathbb E^{x,v}[f(x_t,v_t)] - \pi(f),
\end{equation*}
where $\mathbb E^{x,v}$ means the dynamics $(x_t,v_t)_{t\Ge0}$ starts with the initial state $(x,v)$.

Then we employ the coupling contractivity of the underdamped Langevin dynamics \eqref{underdamped} prove $|\nabla u(x,v,t)|$ have exponential decay in time (see Lemma~\ref{lemma: estimate u 1}).
\item Utilize the function $u(x,v,t)$ to construct the discrete Poisson solution:
\begin{equation*}
	\phi_h(x,v) = h \sum_{n=0}^\infty u(x,v,nh),
\end{equation*}
where $h>0$ represents the time step of the numerical method. Formally, $\phi_h(x,v)$ serves as the solution to the discrete Poisson equation:
\begin{equation*}
	\frac{1-e^{h\mathcal L}}{h} \phi_h(x,v) = f(x,v) - \pi(f),
\end{equation*}
where $\mathcal L$ denotes the generator of \eqref{underdamped}. With the estimates of 
 $|u(x,v,t)|$ and $|\nabla u(x,v,t)|$, $|\nabla \phi_h(x,v)|$ can be globally bounded across the variables $x,v$.
\item Consider the numerical solution $(X_n,V_n)_{n\Ge0}$ with time step $h$, and let $r\Ge1$. Utilizing the Lyapunov function $\mathcal H(x,v)$ in \eqref{Lyapunov}, the uniform-in-time moment estimate
\begin{equation*}
	\sup_{n\Ge0} \big(|X_n| + |V_n| + 1\big)^{2r} \Le C
	\big(|X_0| + |V_0| + 1\big)^{2r}
\end{equation*}
holds true for some constant $C$ independent of the time step $h$.
\item For the numerical solution $Z_n = (X_n,V_n)$, the difference $\phi_h(Z_{n+1}) - \phi_h(Z_n)$ satisfies
\begin{equation}
\begin{aligned}
& \frac{\phi_h(Z_{n+1}) - \phi_h(Z_n)}{h} + f(Z_n) - \pi(f) = \\
& \underbrace{\frac{\phi_h(Z_{n+1}) - \phi_h(Z_n(h))}{h}}_{S_n:\text{ local error}}
+ 
\underbrace{\frac{\phi_h(Z_n(h)) - \phi_h(Z_n)}{h} + f(Z_n) - \pi(f)}_{T_n:\text{ mean-zero}}.
\end{aligned}
\label{rountine}
\end{equation}
Here, the term $S_n$ corresponds to the local error of the numerical integrator, and the term $T_n$ is a mean-zero random variable, because
\begin{equation*}
	\mathbb E\big[T_n|Z_n\big] = 
	\mathbb E\bigg[
	\frac{(e^{h\mathcal L}\phi_h)(Z_n) - \phi_h(Z_n)}h
	\bigg] + f(Z_n) - \pi(f)
	 = 0.
\end{equation*}
As a consequence, the random variables $\{T_n\}_{n=0}^\infty$ are mutually independent, i.e., $\mathbb E[T_n T_m] = 0$ for distinct $n,m\Ge 0$.
\item The summation of \eqref{rountine} over $n$ gives
\begin{equation}
	e(N,h) = \frac1{N}\sum_{n=0}^{N-1}
	 f(Z_n) - \pi(f) = 
	 \frac{\phi_h(Z_0) - \phi_h(Z_N)}{Nh} + 
	 \frac1N \sum_{n=0}^{N-1} (S_n + T_n).
	 \label{rountine 2}
\end{equation}
Then using Cauchy's inequality,
\begin{align}
	\mathbb E[e^2(N,h)] &\Le  
	3\mathbb E\Bigg[\bigg|
	\frac{\phi_h(Z_0) - \phi_h(Z_N)}{Nh}
	\bigg|^2\Bigg] + 3 \mathbb E\Bigg[
	\bigg|\frac1N \sum_{n=0}^{N-1} S_n \bigg|^2
	\Bigg] + 3 \mathbb E\Bigg[
	\bigg|\frac1N \sum_{n=0}^{N-1} T_n \bigg|^2
	\Bigg] \notag \\
	& \Le
	\frac{3}{N^2h^2}
	\mathbb E\big[|\phi_h(Z_0) - \phi_h(Z_N)|^2\big] + \frac3N 
	\sum_{n=0}^{N-1} \mathbb E[S_n^2] + 
	\frac3{N^2} \sum_{n=0}^{N-1} \mathbb E[T_n^2],
	\label{rountine 3}
\end{align}
yielding the quantitative estimate of the statistical error $\mathbb E[e^2(N,h)]$.
\end{enumerate}
When the evolution time $Nh\gg1$, $\mathbb{E}[e^2(N,h)]$ primarily stems from the time discretization error in $\mathbb{E}[S_n^2]$, which necessitates the estimation of the discrete Poisson solution $\phi_h(x,v)$. When an estimate of $|\nabla \phi_h(x,v)|$ is available, an $\mathcal{O}(h^{2p-1})$ error, as stated in Theorem~\ref{theorem: numerical error}, can be achieved. When an estimate of $|\nabla^2 \phi_h(x,v)|$ is attainable, an $\mathcal{O}(h^{2p})$ error, as stated in Theorem~\ref{theorem: numerical error convex}, can be attained.

Two significant advancements characterize our numerical analysis methodology. Firstly, we establish the boundedness of the second derivative $\nabla^2 u(x,v,t)$ when the potential function $U(x)$ exhibits global convexity (refer to Lemma~\ref{lemma: u2 estimate}). This property stems from a thorough analysis of the contractivity in the tangent process, as delineated in Lemmas~\ref{lemma: qp decay} and \ref{lemma: qp diff decay}. Secondly, in the error analysis framework, we employ the discrete Poisson solution $\phi_h(x,v)$ rather than the infinite integral $\phi(x,v) = \int_0^\infty u(x,v,t) \d t$ as commonly used in prior works \cite{Poisson_1, Poisson_2}, to assess the statistical error. Notably, the discrete Poisson solution $\phi_h(x,v)$ ensures that $T_n$ is a mean-zero random variable, thereby rendering $\{T_n\}_{n=0}^{N-1}$ are mutually independent random variables.

\section{Numerical integrator and stability}
\label{section: 3}
We introduce the definition of the strong order of a numerical integrator for underdamped Langevin dynamics \eqref{underdamped}.
For convenience, in Sections~\ref{section: 3} to \ref{section: 5}, we exclusively investigate numerical integrators driven by the deterministic gradient $\nabla U(x)$. The scenario involving the stochastic gradient $b(x,\omega)$ will be addressed comprehensively in Section~\ref{section: 6}.

The strong order of a numerical integrator is defined as follows.
\begin{definition}
\label{defi 1}
Suppose $p>0$ and $q\Ge 1$. Given the intial state $Z_0 = (X_0,V_0)$ in $\mathbb R^{2d}$, the numerical integrator
$$
	Z_0 = (X_0,V_0) \mapsto Z_1 = (X_1,V_1)
$$
has \emph{strong order $(p,q)$} if
\begin{itemize}
\item For some constant $C>0$, the local error can be decomposed as
\begin{equation}
	Z_1 - Z_0(h) = M_h + N_h,
	\label{strong error form}
\end{equation}
and the random variables $M_h$ and $N_h$ satisfy
\begin{equation}
	\mathbb E[M_h] = 0,~~~
	\mathbb E\big[|M_h|^4\big] \Le Ch^{4p+2}(|Z_0|+1)^{4q},~~~
	\mathbb E\big[|N_h|^4\big] \Le Ch^{4p+4}(|Z_0|+1)^{4q}.
	\label{cond 1}
\end{equation}
\item For some constant $C>0$, the local error is bounded by
\begin{equation}
	\mathbb E|Z_1 - Z_0(h)|^4 \Le Ch^6(|Z_0|+1)^4.
	\label{cond 2}
\end{equation}
\item For any $r\Ge 1$, there exists a constant $C>0$ depending on $r$ such that
\begin{equation}
	\mathbb E|Z_1|^r \Le C(|Z_0|+1)^r.
	\label{cond 3}
\end{equation}
\end{itemize}
Here, $p>0$ is the usual strong order of the numerical integrator, and $q\Ge1$ signifies the order of moments required on the initial state $Z_0 = (X_0,V_0)$.
A numerical integrator has \emph{strong order $p$} if there exists $q\Ge 1$ such that this integrator has strong order $(p,q)$.
\end{definition}
\begin{remark}
Remarks on the definition of the strong order:
\begin{enumerate}
\item In the local error $Z_1 - Z_0(h) = M_h + N_h$,  $M_h$ represents the mean-zero variable, while $N_h$ denotes the high-order term. Usually, $M_h$ constitutes a stochastic integral with respect to the Brownian motion $B_t$ (see the proof of Lemma~\ref{lemma: strong order} in Appendix), consequently manifesting as a normally distributed random variable.
\item The condition \eqref{cond 1} is a stronger requirement compared to the conventional definition of strong order, as elucidated in Theorem~1.1 of \cite{Milstein} and Assumption~22 of \cite{UBU_2}. This heightened rigor arises from our requirement of the fourth moment constraint on the local error, as opposed to the second moment typically required. Indeed, the second moment bound
\begin{equation*}
\mathbb E\big[|M_h|^2\big] \Le C h^{2p+1} (|Z_0|+1)^{2q},
~~~~
\mathbb E\big[|N_h|^2\big] \Le C h^{2p+2} (|Z_0|+1)^{2q},
\label{local 1}
\end{equation*}
as delineated in Equation (22) of \cite{UBU_2}, can be can be readily derived from \eqref{cond 1}.
\item The conditions \eqref{cond 2} and \eqref{cond 3} ensure the consistency and the linear stability of the numerical integrator (see Condition~7.1 of \cite{Lyapunov}). Moreover, when $p=q=1$, \eqref{cond 1} directly implies \eqref{cond 2}.
\end{enumerate}
\end{remark}

For the underdamped Langevin dynamics \eqref{underdamped}, various numerical integrators have been developed both for practical simulation and theoretical analysis. These include the Euler--Maruyama integrator, the left-point integrator \cite{left_1,left_2,left_3}, the family of UBU integrators \cite{UBU_1,UBU_2,UBU_3}, and the family of BAOAB integrators \cite{BAOAB_1,BAOAB_2,BAOAB_3}.
In this paper, we demonstrate the application of the error analysis framework by the Euler--Maruyama integrator \eqref{EM} and the UBU integrator \eqref{UBU}, which necessitate only a single evaluation of the gradient at each iteration.
\begin{definition}
Consider the following numerical integrators for \eqref{underdamped}.
\begin{enumerate}
\item \emph{Euler--Maruyama integrator}.
\begin{equation}
\left\{
\begin{aligned}
X_{n+1} & = X_n + hV_n, \\
V_{n+1} & = V_n - h \nabla U(X_n) - \gamma h V_n + \sqrt{2\gamma} \xi_n,
\end{aligned}
\right.
\tag{EM}
\label{EM}
\end{equation}
where the random variable $\xi_n\sim\mathcal N(0,I_d)$ is sampled independently at each time step.
\item \emph{UBU integrator}. We decompose \eqref{underdamped} into two dynamics:
\begin{equation*}
	\mathcal U:\left\{
	\begin{aligned}
	\dot x_t & = v_t, \\
	\dot v_t & = -\gamma v_t + \sqrt{2\gamma} \dot B_t,
	\end{aligned}
	\right.~~~~
	\mathcal B:\dot v_t = -\nabla U(x_t),
\end{equation*}
where both the $\mathcal{U}$ and $\mathcal{B}$ dynamics can be explicitly solved (refer to Equations (2.2)-(2.4) of \cite{UBU_3}). The update from $(X_n,V_n)$ to $(X_{n+1},V_{n+1})$ is given by
\begin{equation}
	(X_n,V_n)
	\xrightarrow{\mathcal U\mathrm{\,for\,}\frac h2}
	(Y_n,V_n^*)
	\xrightarrow{\mathcal B\mathrm{\,for\,} h}
	(Y_n,V_n^{**}) 
	\xrightarrow{\mathcal U\mathrm{\,for\,}\frac h2}
	(X_{n+1},V_{n+1}).
\tag{UBU}
\label{UBU}
\end{equation}
\end{enumerate}
\end{definition}
If we couple \eqref{EM} and the underdamped Langevin dynamics \eqref{underdamped} synchronously, the random variable $\xi_n = (B_{(n+1)h} - B_{nh})/\sqrt{h}$, where $(B_t)_{t\Ge0}$ is the Brownian motion.

The strong order of the numerical integrators \eqref{EM} and \eqref{UBU} is stated as follows.
\begin{lemma}
\label{lemma: strong order}
Assume \textbf{A1}. In the sense of Definition~\ref{defi 1},
\begin{enumerate}
\setzero
\item The Euler--Maruyama integrator \eqref{EM} has strong order $(1,1)$.
\item The UBU integrator \eqref{UBU} has strong order $(2,2)$.
\end{enumerate}
\end{lemma}
The proof of Lemma~\ref{lemma: strong order} is provided in Appendix, wherein the demonstration of the strong order of \eqref{UBU} is adapted from Section 7.6 of \cite{UBU_2}.

In addition to the strong order of the local error, the uniform-in-time moment estimate also plays a crucial role in proving the convergence of the numerical integrators.
\begin{lemma}
\label{lemma: moment}
Assume \textbf{A1}, \textbf{A2}. Suppose the numerical integrator has strong order $(p,q)$, and let $r\Ge 1$ be an arbitrary constant. For some constants $\gamma_0,h_0>0$, whenever $\gamma\Ge \gamma_0$, there exists a constant $C$ depending on $r$ and $\gamma$ such that for any $0<h\Le h_0$,
\begin{equation}
	\sup_{n\Ge0} \mathbb E\big(|X_n| + |V_n| + 1\big)^{2r}
	\Le C\mathbb E\big(|X_0| + |V_0| + 1\big)^{2r}.
	\label{uniform moment bound}
\end{equation}
\end{lemma}
Lemma~\ref{lemma: moment} guarantees that the numerical solution $(X_n,V_n)_{n\Ge0}$ consistently exhibits bounded moments, thereby enabling the application of the local strong error estimate \eqref{cond 1} in subsequent error analyses. A pivotal aspect of the proof lies in the establishment of
\begin{equation}
	\mathcal H(x,v) = \begin{pmatrix}
	x^\T & v^\T 
	\end{pmatrix}
	\begin{pmatrix}
	\gamma & 1 \\
	1 & 1
	\end{pmatrix}
	\begin{pmatrix}
	x \\ v
	\end{pmatrix} = \gamma|x|^2 + 2x^\T v + |v|^2 ,
	\label{Lyapunov}
\end{equation}
which serves as a Lyapunov function for the underdamped Langevin dynamics \eqref{underdamped}. 
The construction of $\mathcal{H}(x,v)$ draws inspiration from Theorem 3.5 of \cite{Pavliotis_2} and Theorem 1 of \cite{coupling_1}. The complete proof of Lemma~\ref{lemma: moment} is provided in Appendix.
\begin{remark}
Lemma~\ref{lemma: moment} necessitates the damping rate $\gamma$ to be sufficiently large. For a general $\gamma>0$, an alternative choice for the Lyapunov function is available, such as
\begin{equation*}
\mathcal H(x,v) = \frac12|v|^2 + U(x) + \frac{\gamma}2
x^\T v + \frac{\gamma^2}4|x|^2 + 1
\end{equation*}
as proposed in \cite{Lyapunov,coupling_0}. Nevertheless, we favor \eqref{Lyapunov} due to its exclusion of the potential function $U(x)$, rendering the error analysis methodology suitable for numerical integrators featuring stochastic gradients.
\end{remark}
A corollary of Lemma~\ref{lemma: moment} yields the following diffusion estimate, which proves useful in the demonstration of Theorems~\ref{theorem: numerical error} and \ref{theorem: numerical error convex}.
\begin{lemma}
\label{lemma: diffusion}
Assume \textbf{A1}, \textbf{A2}. Suppose the numerical integrator has strong order $(p,q)$. For some constants $\gamma_0,h_0>0$, whenever $\gamma\Ge \gamma_0$, there exists a constant $C$ depending on $\gamma$ such that for any $0<h\Le h_0$ and integer $N\Ge0$,
\begin{equation}
	\mathbb E\big[|Z_N - Z_0|^2\big] \Le CNh \mathbb E
	(|Z_0|+1)^{2q}.
\end{equation}
\end{lemma}
In the proof of Lemma~\ref{lemma: diffusion}, we decompose $Z_1 - Z_0$ into the summation of a mean-zero random variable $\bar M_h$ and a high-order term $\bar N_h$, satisfying
\begin{equation*}
\mathbb E\big[|\bar M_h|^2\big] \Le C h (|Z_0|+1)^{2q},~~~~
\mathbb E\big[|\bar N_h|^2\big] \Le Ch^2 
(|Z_0|+1)^{2q}.
\end{equation*}
The detailed proof of Lemma~\ref{lemma: diffusion} is provided in Appendix.

\section{Estimate of the Kolmogorov equation}
\label{section: 4}
In this section, we delve into the pointwise estimate of the solution $u(x,v,t)$ to the Kolmogorov equation, which holds paramount importance in constructing the discrete Poisson solution outlined in Section~\ref{section: 5}. Leveraging the Wasserstein-1 contractivity of \eqref{underdamped}, we bound $|u(x,v,t)|$ and $|\nabla u(x,v,t)|$ when $U(x)$ is strongly convex outside a ball. Furthermore, we utilize synchronous coupling techniques to estimate $|\nabla^2 u(x,v,t)|$ when $U(x)$ is strongly convex in $\mathbb{R}^d$.

\subsection{Estimate under strong convexity outside a ball}

We introduce the Wasserstein-1 contractivity property of the underdamped Langevin dynamics \eqref{underdamped}, which holds pivotal significance in estimating the Kolmogorov equation. We denote the Markov semigroup of \eqref{underdamped} by $(p_t)_{t\Ge0}$, wherein for any initial distribution $\nu\in\mathcal{P}(\mathbb{R}^{2d})$, the resultant distribution $\nu p_t\in\mathcal{P}(\mathbb{R}^{2d})$ signifies the distribution law of $(x_t,v_t)$ provided $(x_0,v_0)\sim \nu$. Moreover, the Wasserstein-1 distance in $\mathbb{R}^{2d}$ is 
\begin{equation*}
\mathcal W_1(\mu,\nu) = 
\inf_{\gamma\in\Pi(\mu,\nu)}
\int_{\mathbb R^{2d}}
\big(|x-x'| + |v-v'|\big) 
\gamma(\d x\d v\d x'\d v'),
\end{equation*}
where $\Pi(\mu,\nu)$ is set of the transport plans between the distributions $\mu$ and $\nu$.
\begin{theorem}[contractivity]
\label{theorem: contractivity}
Assume \textbf{A1}, \textbf{A2}. For some constant $\gamma_0 >0$, whenever the damping rate $\gamma\Ge \gamma_0$, there exist constants $C,\lambda>0$ depending on $\gamma$ such that
\begin{equation}
\mathcal W_1(\mu p_t,\nu p_t) \Le Ce^{-\lambda t}
\mathcal W_1(\mu,\nu),~~~~t\Ge 0,
\end{equation}
for any probability distributions $\mu,\nu \in \mathcal P(\mathbb R^{2d})$.
\end{theorem}
The proof of Theorem~\ref{theorem: contractivity} using the reflection coupling can be found in Theorem~5 of \cite{coupling_1}.
\begin{remark}
The first reflection coupling result regarding the underdamped Langevin dynamics \eqref{underdamped} is provided in Theorem 2.3 of \cite{coupling_0}, albeit solely in the context of a semimetric. In comparison to the Wasserstein-1 distance, contractivity in a semimetric poses challenges in establishing estimates for the Kolmogorov equation.
\end{remark}

Next, we present the estimation of the Kolmogorov equation when the potential function $U(x)$ is strongly convex outside a ball. The generator of the underdamped Langevin dynamics \eqref{underdamped} is expressed as
\begin{equation*}
\mathcal L = v\cdot \nabla_x - (\nabla U(x) + \gamma v)
\cdot \nabla_v + \gamma \Delta_v.
\end{equation*}
We define $u(x,v,t)$ to be the solution of the following Kolmogorov equation
\begin{equation*}
\left\{
\begin{aligned}
\frac{\partial u}{\partial t}(x,v,t) & = (\mathcal L u)
(x,v,t), && t\Ge0,\\
u(x,v,0) & = f(x,v) - \pi(f),
\end{aligned}
\right.
\end{equation*}
where $\pi(f)$ represents the average value of $f(x,v)$. Then, $u(x,v,t)$ can be expressed as
\begin{equation}
	u(x,v,t) = (e^{t\mathcal L}f)(x,v) - \pi(f) = 
	\mathbb E^{x,v}[f(x_t,v_t)] - \pi(f),
	\label{u expression}
\end{equation}
where the superscript $x,v$ denotes that the initial state of $(x_t,v_t)_{t\Ge0}$ is $(x,v)\in\mathbb{R}^{2d}$.

The Wasserstein-1 contractivity outlined in Theorem~\ref{theorem: contractivity} offers a straightforward estimation of $u(x,v,t)$ and its derivatives. This approach, known as Stein's method, provides a means of estimating $u(x,v,t)$ utilizing contractivity properties. For further insights, refer to Theorem 5 of \cite{Poisson_3}.
\begin{lemma}
\label{lemma: estimate u 1}
Assume \textbf{A1}, \textbf{A2}, \textbf{A3}. Let $u(x,v,t)$ be the Kolmogorov solution defined in \eqref{u expression}. For some constant $\gamma_0>0$, whenever $\gamma\Ge \gamma_0$, there exist constants $C,\lambda>0$ depending on $\gamma$ such that for any $x,v\in\mathbb R^d$,
\begin{equation}
|u(x,v,t)| \Le C e^{-\lambda t}(|x|+|v|+1),~~~~
|\nabla u(x,v,t)| \Le Ce^{-\lambda t}.
\label{u estimate}
\end{equation}
\end{lemma}
The proof of Lemma~\ref{lemma: estimate u 1} is provided in Appendix. Although $u(x,v,t)$ and $\nabla u(x,v,t)$ can be readily bounded using the Wasserstein-1 contractivity, estimating the second derivative $\nabla^2 u(x,v,t)$ remains challenging without the global convexity of the potential $U(x)$.
\begin{remark}
The availability of the estimate for the second derivative $\nabla^2 u(x,v,t)$ is the distinguishing factor leading to the difference in the order of the time step $h$ between Theorems~\ref{theorem: numerical error} and \ref{theorem: numerical error convex}. It is remarkable that the second derivative estimate in the overdamped Langevin dynamics can be readily accomplished from the $L^p$-estimate for elliptic operators, as elaborated in Lemma 17 of \cite{Poisson_2} and Theorem 5 of \cite{Poisson_3}.
\end{remark}
\subsection{Estimate under strong convexity in $\mathbb R^d$}
We proceed to provide an estimate of $\nabla^2 u(x,v,t)$ under the assumption that the potential function $U(x)$ is strongly convex in $\mathbb{R}^d$. Our approach relies on the tangent process introduced in Section 3.2 of \cite{Pavliotis_2} to characterize the derivative $\nabla u(x,v,t)$. The tangent process, also referred to as the first variation process \cite{variation_1,variation_2,variation_3}, is instrumental in our analysis.
\begin{definition}
Given the initial state $(x,v)$ and the Brownian motion $(B_t)_{t\Ge0}$, let $(x_t,v_t)_{t\Ge0}$ be the strong solution to \eqref{underdamped} with initial state $(x,v)\in\mathbb R^{2d}$. The tangent processes
$$
D_x x_t,D_x v_t,D_v x_t,D_v v_t\in\mathbb R^{d\times d}
$$
are defined as the partial derivatives of $(x_t,v_t)$ with respect to $(x,v)$.
\end{definition}
By differentiating \eqref{underdamped} with respect to $(x,v)$, we derive that $(D_x x_t, D_x v_t)_{t\Ge0}$ satisfies
\begin{equation}
	\frac{\partial}{\partial t}
	\begin{pmatrix}
	D_x x_t \\ D_x v_t
	\end{pmatrix} = 
	\begin{pmatrix}
	O_d & I_d \\
	-\nabla^2 U(x_t) & -\gamma I_d
	\end{pmatrix}
	\begin{pmatrix}
	D_x x_t \\ D_x v_t
	\end{pmatrix},~~~~
	\begin{pmatrix}
	D_x x_0 \\ D_x v_0
	\end{pmatrix} = 
	\begin{pmatrix}
	I_d \\ O_d
	\end{pmatrix},
	\label{Dx process}
\end{equation} 
and $(D_v x_t, D_v v_t)_{t\Ge0}$ satisfies
\begin{equation}
	\frac{\partial}{\partial t}
	\begin{pmatrix}
	D_v x_t \\ D_v v_t
	\end{pmatrix} = 
	\begin{pmatrix}
	O_d & I_d \\
	-\nabla^2 U(x_t) & -\gamma I_d
	\end{pmatrix}
	\begin{pmatrix}
	D_v x_t \\ D_v v_t
	\end{pmatrix},~~~~
	\begin{pmatrix}
	D_v x_0 \\ D_v v_0
	\end{pmatrix} = 
	\begin{pmatrix}
	O_d \\ I_d
	\end{pmatrix}.
	\label{Dv process}
\end{equation}
In a unified framework, the matrix-valued process $(Q_t,P_t)_{t\Ge0}$ defined by
\begin{equation}
\begin{pmatrix}
\dot Q_t \\ \dot P_t
\end{pmatrix} = 
\begin{pmatrix}
O_d & I_d \\ -\nabla^2 U(x_t) & -\gamma I_d
\end{pmatrix}
\begin{pmatrix}
Q_t \\ P_t
\end{pmatrix}
\label{qp process}
\end{equation}
encompasses both tangent processes \eqref{Dx process} and \eqref{Dv process}. It's noteworthy that $(Q_t,P_t)_{t\Ge0}$ solely depends on the initial state $(x,v)$, the Brownian motion $(B_t)_{t\Ge0}$, and the initial value $(Q_0,P_0)$. We can thus represent a tangent process using the following mapping:
\begin{equation}
(Q_t,P_t)_{t\Ge0} = \mathrm{TangentProcess}\big[
(x,v),
(B_t)_{t\Ge0},(Q_0,P_0)
\big].
\label{qp process mapping}
\end{equation}
The tangent process mapping \eqref{qp process mapping} exhibits linearity in $(Q_0,P_0)$, yet it relies nonlinearly on the initial state $(x,v)$. This nonlinearity arises from the fact that the dynamics \eqref{qp process} directly depend on the trajectory $(x_t,v_t)_{t\Ge0}$.

The tangent process provides a straightforward representation of $\nabla u(x,v,t)$ with the Kolmogorov solution $u(x,v,t)$ defined in \eqref{u expression}. By taking the derivative with respect to the initial state $(x,v)$ in \eqref{u expression}, we obtain the following expressions:
\begin{subequations}
\begin{align}
	\nabla_x u(x,v,t) & = \mathbb E^{x,v}
	\Big[
	\nabla_x f(x_t,v_t) \cdot D_x x_t + 
	\nabla_v f(x_t,v_t) \cdot D_x v_t
	\Big], \\
	\nabla_v u(x,v,t) & = \mathbb E^{x,v}
	\Big[
	\nabla_x f(x_t,v_t) \cdot D_v x_t + 
	\nabla_v f(x_t,v_t) \cdot D_v v_t
	\Big].
\end{align}
\label{grad u expression}
\end{subequations}
The expressions \eqref{grad u expression} reveal that the exponential decay of $\nabla u(x,v,t)$ can be inferred from the exponential decay of the tangent processes $D_x x_t,D_x v_t,D_v x_t,D_v v_t$.
\begin{lemma}
\label{lemma: qp decay}
Assume \textbf{A1}, \textbf{A2'}, \textbf{A3}. For some constant $\gamma_0>0$, whenever $\gamma\Ge \gamma_0$, there exist constants $C,\lambda>0$ depending on $\gamma$ such that for any initial state $(x,v)$, the Brownian motion $(B_t)_{t\Ge0}$ and the initial value $(Q_0,P_0)$, the tangent process $(Q_t,P_t)_{t\Ge0}$ defined by
\begin{equation*}
(Q_t,P_t)_{t\Ge0} = \mathrm{TangentProcess}\big[
(x,v),
(B_t)_{t\Ge0},(Q_0,P_0)
\big]
\end{equation*}
satisfies
\begin{equation}
|Q_t|+|P_t| \Le C e^{-\lambda t} (|Q_0|+|P_0|).
\end{equation}
\end{lemma}
Inspired from the proof of Lemma~\ref{lemma: moment}, we consider the Lyapunov function
\begin{equation*}
\mathcal H(Q,P) := 
\Tr[
\begin{pmatrix}
Q^\T & P^\T
\end{pmatrix}
	\begin{pmatrix}
	\gamma I_d & I_d \\
	I_d & I_d
	\end{pmatrix}
\begin{pmatrix}
Q \\ P
\end{pmatrix}] = 
\Tr[\gamma Q^\T Q + 2 Q^\T P + P^\T P],
\end{equation*}
and prove Lemma~\ref{lemma: qp decay} via the exponential decay of $\mathcal H(Q_t,P_t)$.
The complete proof of Lemma~\ref{lemma: qp decay} is provided in Appendix.

The exponential decay of the tangent process $(Q_t,P_t)_{t\Ge0}$ is adequate for estimating $\nabla u(x,v,t)$ as given in \eqref{grad u expression}. However, to estimate $\nabla^2 u(x,v,t)$, more comprehensive results on the tangent process $(Q_t,P_t)$ are necessary. Specifically, we require an analysis of the coupling error between two tangent processes driven by different initial states.
\begin{lemma}
\label{lemma: qp diff decay}
Assume \textbf{A1}, 
\textbf{A2'}, \textbf{A3}. For some constant $\gamma_0>0$, whenever $\gamma\Ge \gamma_0$, there exist constants $C,\lambda>0$ depending on $\gamma$ such that 
for any initial states $(x,v)$ and $(x',v')$, the Brownian motion $(B_t)_{t\Ge0}$ and the initial value $(Q_0,P_0)$, the tangent processes $(Q_t,P_t)_{t\Ge0}$ and $(Q_t',P_t')_{t\Ge0}$ defined by
\begin{align*}
(Q_t,P_t)_{t\Ge0} & = \mathrm{TangentProcess}\big[
(x,v),
(B_t)_{t\Ge0},(Q_0,P_0)
\big], \\
(Q_t',P_t')_{t\Ge0} & = \mathrm{TangentProcess}\big[
(x',v'),
(B_t)_{t\Ge0},(Q_0,P_0)
\big],
\end{align*}
satisfy
\begin{equation}
|Q_t - Q_t'| + |P_t - P_t'| \Le C e^{-\lambda t}
(|Q_0| + |P_0|) \big(
|x-x'| + |v-v'|
\big).
\label{qp diff decay}
\end{equation}
\end{lemma}
The proof of Lemma~\ref{lemma: qp diff decay} is provided in the Appendix and relies on the exponential decay of the Lyapunov function $\mathcal H(Q_t - Q_t',P_t - P_t')$.

Formally, Lemma~\ref{lemma: qp diff decay} implies the derivatives of $D_xx_t,D_xv_t,D_vx_t,D_vv_t$ have exponential decay, and thus we can the derivative with respect to $(x,v)$ in $\nabla u(x,v,t)$ in \eqref{grad u expression} to show that $\nabla^2 u(x,v,t)$ has exponential decay.
\begin{lemma}
\label{lemma: u2 estimate}
Assume \textbf{A1}, \textbf{A2'}, \textbf{A3}.
Let $u(x,v,t)$ be the Kolmogorov solution defined in \eqref{u expression}.
For some constant $\gamma_0>0$, whenever $\gamma\Ge \gamma_0$, there exist constants $C,\lambda>0$ depending on $\gamma$ such that for any $x,v\in\mathbb R^d$,
\begin{equation}
	|\nabla^2 u(x,v,t)| \Le Ce^{-\lambda t}.
	\label{u2 estimate}
\end{equation}
\end{lemma}
The proof of Lemma~\ref{lemma: u2 estimate} is given in Appendix.

When $U(x)$ lacks global convexity, obtaining estimates for $\nabla^2 u(x,v,t)$ becomes challenging. The estimation of $\nabla^2 u(x,v,t)$ in cases of global convexity heavily depends on the exponential decay of the tangent process $(Q_t,P_t)$, a property whose verification remains intractable when $U(x)$ is nonconvex.
\section{Analysis of statistical error with deterministic gradients}
\label{section: 5}
In this section, we analyze the statistical error of the numerical integrators outlined in Section~\ref{section: 3}, focusing on scenarios where the drift force is governed by the deterministic gradient $\nabla U(x)$. Our approach begins with leveraging the Kolmogorov solution $u(x,v,t)$ discussed in Section~\ref{section: 4} to construct the discrete Poisson solution $\phi_h(x,v,t)$. Subsequently, we derive representations for the time average of the numerical solution $(X_n,V_n)_{n\Ge0}$.
\subsection{Construction of the discrete Poisson solution}
The discrete Poisson solution $\phi_h(x,v)$ is defined by the expression:
\begin{equation}
\phi_h(x,v) = h\sum_{n=0}^\infty u(x,v,nh),
\label{phi_h expression}
\end{equation}
where $u(x,v,t)$ represents the Kolmogorov solution in \eqref{u expression}. Given that $u(x,v,t) = (e^{t\mathcal L} f)(x,v) - \pi(f)$, the function $\phi_h(x,v)$ can be formally expressed as
\begin{equation*}
	\phi_h(x,v) = h\sum_{n=0}^\infty e^{nh\mathcal L}\big(f(x,v) - \pi(f)\big) = 
	h\big(1 - e^{h\mathcal L}\big)^{-1} f,
\end{equation*}
leading to the derivation of the discrete Poisson equation,
\begin{equation}
	\frac{1-e^{h\mathcal L}}h \phi_h(x,v) = f(x,v) - \pi(f).
\label{discrete Poisson}
\end{equation}
The term \emph{discrete} is employed here to underscore the distinction from the integral
\begin{equation*}
\phi(x,v) = \int_0^\infty u(x,v,t),
\end{equation*}
which represents the solution to the Poisson equation
$$
	-\mathcal L \phi(x,v) = f(x,v) - \pi(f).
$$
Evidently, as the time step $h\rightarrow 0$, $\phi(x,v)$ emerges as the formal limit of $\phi_h(x,v)$.

The Poisson solution $\phi(x,v)$ has found widespread application in scrutinizing the statistical error across various explicit numerical integrators \cite{Poisson_1,Poisson_2,Poisson_3,Poisson_4}, in particular the (stochastic gradient) Euler--Maruyama integrator for the overdamped Langevin dynamics. Despite its efficacy, employing the Poisson equation approach often necessitates high-order derivatives of the solution $\phi(x,v)$ (refer, for instance, to Section~5.2 of \cite{SGLD_3}) and mandates acquaintance with the specific form of the numerical integrator.

Conversely, in this paper, we adhere to the approach outlined in Section~\ref{section: 2} and leverage the discrete Poisson solution $\phi_h(x,v)$ rather than $\phi(x,v)$ in analyzing the statistical error. This facilitates a more straightforward investigation of a broad spectrum of numerical integrators characterized by a time step $h>0$.

Given the exponential decay of the function $u(x,v,t)$ in Lemmas~\ref{lemma: estimate u 1} and \ref{lemma: u2 estimate}, we can readily establish the estimates of the discrete Poisson solution $\phi_h(x,v)$ defined in \eqref{phi_h expression} in a convenient manner.
\begin{lemma}
\label{lemma: phi h 1} Assume \textbf{A1}, \textbf{A2}, \textbf{A3}.
Let $\phi_h(x,v)$ be the discrete Poisson solution defined in \eqref{phi_h expression}.
For some constant $\gamma_0$, whenever $\gamma\Ge \gamma_0$, there exists a constant $C$ depending on $\gamma$ such that for any $0<h\Le 1$ and $x,v\in\mathbb R^d$,
\begin{equation}
	|\nabla \phi_h(x,v)| \Le C.
\end{equation}
\end{lemma}
\begin{lemma}
\label{lemma: phi h 2}
Assume \textbf{A1}, \textbf{A2'}, \textbf{A3}. 
Let $\phi_h(x,v)$ be the discrete Poisson solution defined in \eqref{phi_h expression}.
For some constant $\gamma_0$, whenever $\gamma\Ge \gamma_0$, there exists a constant $C$ depending on $\gamma$ such that for any $0<h\Le 1$ and $x,v\in\mathbb R^d$,
\begin{equation}
	|\nabla^2\phi_h(x,v)| \Le C.
\end{equation}
\end{lemma}
Here, the upper bound of the time step $h$ can be replaced by any positive constant. Therefore, there is no inherent constraint on the time step $h$ in Lemmas~\ref{lemma: phi h 1} and \ref{lemma: phi h 2}.
\subsection{Statistical error of numerical integrators}
We examine the numerical solution $(X_n,V_n)_{n\Ge0}$ initialized with $(X_0,V_0)$, employing a numerical integrator with the strong order $(p,q)$. Under \textbf{A2}, we attain the following estimate of $e(N,h)$ as delineated in \eqref{time average error}.
\begin{theorem}
\label{theorem: numerical error}
Assume \textbf{A1}, \textbf{A2}, \textbf{A3}. Suppose the numerical integrator for \eqref{underdamped} has strong order $(p,q)$. For some constant $\gamma_0>0$, whenever $\gamma\Ge \gamma_0$, there exist constants $C,h_0>0$ depending on $\gamma$ such that for any $0<h\Le h_0$,
\begin{equation}
\mathbb E[e^2(N,h)] \Le C\bigg(h^{2p-1}+\frac1{Nh}\bigg)
\mathbb E\big(|X_0|+|V_0|+1\big)^{2q}.
\end{equation}
\end{theorem}
\emph{Proof of Theorem~\ref{theorem: numerical error}. Statistical error of the numerical integrator.}\\[6pt]
Define the quantities $S_n$ and $T_n$ as in \eqref{rountine}, namely,
\begin{equation*}
	S_n = \frac{\phi_h(Z_{n+1}) - \phi_h(Z_n(h))}{h},~~~~
	T_n = \frac{\phi_h(Z_n(h)) - \phi_h(Z_n)}h + 
	f(Z_n) - \pi(f).
\end{equation*}
Next we estimate the error terms $\mathbb E[S_n^2]$ and $\mathbb E[T_n^2]$ respectively.

Using the boundedness of $|\nabla \phi_h(x,v)|$ in Lemma~\ref{lemma: estimate u 1}, $S_n$ is bounded by
\begin{equation*}
	|S_n| \Le \frac{C}h |Z_{n+1} - Z_n(h)|.
\end{equation*}
Since the numerical integrator has strong order $(p,q)$, we obtain
\begin{align}
	\mathbb E[S_n^2] & \Le 
	\frac{C}{h^2} \mathbb E|Z_{n+1} - Z_n(h)|^2 
	\Le 
	Ch^{2p-1}\mathbb E(|Z_n|+1)^{2q} \notag \\
	& \Le Ch^{2p-1}\mathbb E(|Z_0|+1)^{2q}.
	\label{S_n bound 2}
\end{align}
In the last inequality, we have used the uniform-in-time $2q$-th moment estimate in Lemma~\ref{lemma: moment}.

Again using the boundedness of $|\nabla \phi_h(x,v)|$, $T_n$ is bounded by
\begin{equation*}
	|T_n| \Le \frac{C|Z_n(h) - Z_n|}{h} + 
	|f(Z_n)| + C.
\end{equation*}
When the time step $h\Le h_0$, from It\^o's calculus we have the inequality
$$
\mathbb E|Z_n(h) - Z_n|^2 \Le Ch
\mathbb E(|Z_n|+1)^2 \Le Ch \mathbb E(|Z_0|+1)^2 
,
$$
hence we obtain
\begin{equation}
\mathbb E[T_n^2] \Le \frac{C}{h}\mathbb E(|Z_0|+1)^2.
\label{T_n bound 2}
\end{equation}

From Lemma~\ref{lemma: diffusion}, there holds diffusion estimate of the numerical solution
\begin{equation}
	\mathbb E|Z_N - Z_0|^2 \Le CNh \mathbb E(|Z_0|+1)^{2q}.
	\label{diffusion bound}
\end{equation}
In the representation of the statistical error \eqref{rountine 3}, combining \eqref{S_n bound 2}, \eqref{T_n bound 2} and \eqref{diffusion bound} we obtain
\begin{align*}
\mathbb E[e^2(N,h)] & \Le
\frac{3 C}{N^2h^2}
	\mathbb E\big[|Z_0 - Z_N|^2\big] + \frac3N 
	\sum_{n=0}^{N-1} \mathbb E[S_n^2] + 
	\frac3{N^2} \sum_{n=0}^{N-1} \mathbb E[T_n^2] \\
& \Le \frac{C}{Nh}\mathbb E(|Z_0|+1)^{2q} +  C h^{2p-1}\mathbb E(|Z_0|+1)^{2q} + \frac{C}{Nh} \mathbb E(|Z_0|+1)^2 \notag \\
& \Le C\bigg(h^{2p-1} + \frac1{Nh}\bigg) \mathbb E(|Z_0|+1)^{2q},
\end{align*}
thereby concluding the proof of Theorem~\ref{theorem: numerical error}. \hfill $\blacksquare$\\[6pt]
In the proof of Theorem~\ref{theorem: numerical error}, the estimate of $S_n$ is suboptimal due to the availability of only the first derivative $\nabla \phi_h(x,v)$, yielding an error order of $\mathcal{O}(h^{2p-1})$. To achieve a discretization error of $\mathcal{O}(h^{2p})$, the global convexity of $U(x)$ is required to estimate the second derivative $\nabla^2 \phi_h(x,v)$.
\begin{theorem}
\label{theorem: numerical error convex}
Assume \textbf{A1}, \textbf{A2'}, \textbf{A3}. Suppose the numerical integrator for \eqref{underdamped} has strong order $(p,q)$. For some constant $\gamma_0>0$, whenever $\gamma\Ge \gamma_0$, there exist constants $C,h_0>0$ depending on $\gamma$ such that for any $0<h\Le h_0$,
\begin{equation}
\mathbb E[e^2(N,h)] \Le C\bigg(h^{2p}+\frac1{Nh}\bigg)
\mathbb E\big(|X_0|+|V_0|+1\big)^{4q}.
\end{equation}
\end{theorem}
\emph{Proof of Theorem~\ref{theorem: numerical error convex}. Statistical error of the numerical integrator under global convexity.}\\[6pt]
The proof closely resembles that of Theorem~\ref{theorem: numerical error}; however, a more accurate estimation of the summation $\sum_{n=0}^{N-1} S_n$ is required. We express $S_n$ for each $n\Ge0$ as follows,
\begin{equation}
	S_n = \frac1h \big(Z_{n+1} - Z_n(h)\big) \cdot 
	\int_0^1 \nabla \phi_h\big(
	\theta Z_{n+1} + (1-\theta) Z_n(h)
	\big) \d\theta.
	\label{S_n expression}
\end{equation}
According to Definition~\ref{defi 1}, the local error $Z_{n+1} - Z_n(h) = M_{n,h} + N_{n,h}$ satisfies
\begin{equation*}
\mathbb E\big[M_{n,h}\big|Z_n\big] = 0,~~
\mathbb E\big[|M_{n,h}|^4\big|Z_n\big] \Le C 
h^{4p+2} (|Z_n|+1)^{4q},~~
\mathbb E\big[|N_{n,h}|^4\big|Z_n\big] \Le C
h^{4p+4} (|Z_n|+1)^{4q}.
\end{equation*}
Here, $\mathbb E[\cdot|Z_n]$ means the conditional expectation for given $Z_n = (X_n,V_n)$.
As a consequence, $S_n$ in \eqref{S_n expression} can be decomposed into $S_{n,1} + S_{n,2} + S_{n,3}$, where
\begin{align*}
S_{n,1} & = \frac1h M_{n,h} \cdot \nabla \phi_h(Z_n), \\
S_{n,2} & = \frac1h N_{n,h} \cdot \nabla \phi_h(Z_n), \\
S_{n,3} & = \frac1h (Z_{n+1} - Z_n(h))\cdot 
\int_0^1 \Big(\nabla\phi_h\big(
\theta Z_{n+1} + (1-\theta) Z_n(h)
\big) - \nabla\phi_h(Z_n)\Big)\d\theta.
\end{align*}
Next we estimate the mean square error of $\sum_{n=0}^{N-1} S_{n,i}$ for $i=1,2,3$, respectively.
\begin{enumerate}
\item \emph{Estiamte for $S_{n,1}$.}

Since $\nabla\phi_h(Z_n)$ is globally bounded by Lemma~\ref{lemma: estimate u 1}, we obtain $|S_{n,1}| \Le \frac Ch|M_{n,h}|$. 
Then uniform-in-time moment estimate from Lemma~\ref{lemma: moment} yields
\begin{equation*}
\mathbb E [S_{n,1}^2] \Le 
\frac{C}{h^2} \mathbb E\big[|M_{n,h}|^2\big] \Le
C h^{2p-1}
\mathbb E(|Z_n|+1)^{2q} \Le C h^{2p-1} \mathbb E(|Z_0|+1)^{2q}.
\end{equation*}
Given that $\{S_{n,1}\}_{n=0}^{N-1}$ are unbiased and mutually independent, we arrive at
\begin{equation}
\mathbb E\Bigg[\bigg|\sum_{n=0}^{N-1} S_{n,1}\bigg|^2\Bigg] = \sum_{n=0}^{N-1} \mathbb E [S_{n,1}^2] \Le 
CN h^{2p-1} \mathbb E(|Z_0|+1)^{2q}.
\label{e sum 1}
\end{equation}
\item \emph{Estimate for $S_{n,2}$.}

Again, exploiting the boundedness of $\nabla \phi_h(Z_n)$, we have $|S_{n,2}| \Le \frac{C}{h} |N_{n,h}|$ and
\begin{equation*}
	\mathbb E[S_{n,2}^2] 
	\Le Ch^{2p} \mathbb E(|Z_n|+1)^{2q}
	\Le C h^{2p} \mathbb E(|Z_0|+1)^{2q}.
\end{equation*}
Utilizing Cauchy's inequality, we deduce
\begin{equation}
\mathbb E\Bigg[\bigg|\sum_{n=0}^{N-1} S_{n,2}\bigg|^2\Bigg] \Le N\sum_{n=0}^{N-1} \mathbb E [S_{n,2}^2] \Le 
CN^2 h^{2p} \mathbb E(|Z_0|+1)^{2q}.
\label{e sum 2}
\end{equation}
\item \emph{Estimate for $S_{n,3}$.}

Utilizing the boundedness of $\nabla^2 \phi_h(z)$ from Lemma~\ref{lemma: u2 estimate}, we have
\begin{equation*}
|S_{n,3}| \Le \frac Ch |Z_{n+1} - Z_n(h)|
\big( |Z_{n+1} - Z_n| + |Z_n(h) - Z_n|\big)
\end{equation*}
Applying Cauchy's inequality, we obtain
\begin{equation}
\mathbb E[S_{n,3}^2] \Le \frac C{h^2}
\sqrt{\mathbb E\big[|Z_{n+1} - Z_n(h)|^4\big]
}
\sqrt{\mathbb E\big[|Z_{n+1} - Z_n|^4 + |Z_n(h) - Z_n|^4\big]}
\label{S n3 estimate}
\end{equation}
Utilizing the linear stability condition \eqref{cond 2} in Definition~\ref{defi 1}, we deduce
\begin{align*}
& ~~~~ \mathbb E\big[|Z_{n+1} - Z_n|^4 + |Z_n(h) - Z_n|^4\big]\\
& \Le C\, \mathbb E\big[|Z_{n+1} - Z_n(h)|^4\big] +  C\, \mathbb E\big[|Z_n(h) - Z_n|^4\big]\\
& \Le Ch^6 \mathbb E(|Z_n|+1)^4 + Ch^2\mathbb E(|Z_n|+1)^4 \Le Ch^2 \mathbb E(|Z_0|+1)^4.
\end{align*}
Moreover, the local error bound \eqref{cond 1} implies
\begin{equation*}
\mathbb E\big[|Z_{n+1} - Z_n(h)|^4\big] \Le C h^{4p+2}
\mathbb E(|Z_n|+1)^{4q} \Le C h^{4p+2}
\mathbb E(|Z_0|+1)^{4q}.
\end{equation*}
Therefore, from \eqref{S n3 estimate} and Cauchy's inequality, we conclude
\begin{equation}
\mathbb E\Bigg[\bigg|\sum_{n=0}^{N-1} S_{n,3}\bigg|^2\Bigg]
\Le 
N\sum_{n=0}^{N-1}\mathbb E[S_{n,3}^2] \Le 
CN^2h^{2p} \mathbb E(|Z_0|+1)^{4q}.
\label{e sum 3}
\end{equation}
\end{enumerate}
From the estimates \eqref{e sum 1}, \eqref{e sum 2} and \eqref{e sum 3}, we obtain
\begin{equation*}
\mathbb E\Bigg[\bigg|\sum_{n=0}^{N-1} S_{n}\bigg|^2\Bigg]
\Le 3\sum_{i=1}^3 \mathbb E\Bigg[\bigg|\sum_{n=0}^{N-1} S_{n,i}\bigg|^2\Bigg] \Le 
C(N^2 h^{2p} + N h^{2p-1})\mathbb E(|Z_0|+1)^{4q}.
\end{equation*}
Recalling the expression of the time average error $e(N,h)$ in \eqref{rountine 2}, we obtain
\begin{align}
\mathbb E[e^2(N,h)] & \Le \frac{C}{(Nh)^2}
\mathbb E\big[|Z_N-Z_0|^2\big] + \frac{C}{N^2}
\sum_{n=0}^{N-1}\mathbb E[T_n^2] + \frac{C}{N^2}\mathbb E\Bigg[\bigg|\sum_{n=0}^{N-1} S_{n}\bigg|^2\Bigg]\notag  \\
& \Le \frac{C}{Nh} \mathbb E(|Z_0|+1)^{2q}  + C\bigg( h^{2p} + \frac1{Nh}\bigg)\mathbb E(|Z_0|+1)^{4q} \notag \\
& \Le C\bigg(h^{2p} + \frac1{Nh}\bigg) \mathbb E(|Z_0|+1)^{4q},
\end{align}
which completes the proof of Theorem~\ref{theorem: numerical error convex}. \hfill $\blacksquare$
\section{Analysis of statistical error with stochastic gradients}
\label{section: 6}
In this section, we investigate the statistical error of numerical integrators with stochastic gradients, where the potential function $U(x)$ is strongly convex in $\mathbb R^d$. The stochastic gradient $b(x,\omega)$ involves a random variable $\omega\in\Omega$, which is independently sampled at each time step. 
For example, when $U(x)$ is the average of a set of potential functions $\{U_j(x)\}_{j=1}^M$, namely $U(x) = \frac1M \sum_{j=1}^M U_j(x)$,
a typical choice of the stochastic gradient $b(x,\omega)$ is the mini-batch approximation
\begin{equation}
	b(x,\omega) = \frac1{\mathcal B} \sum_{j\in\mathcal C} 
	\nabla U_j(x),
	\label{mini batch}
\end{equation}
where $\mathcal B\in\mathbb N$ is the batch size, and $\mathcal C = \mathcal C(\omega)$ is a random subset of $\{1,\cdots,M\}$ of size $\mathcal B$. In this paper we adopt the original notation $b(x,\omega)$ rather than the mini-batch approximation \eqref{mini batch} to represent the stochastic gradient.\\[6pt] 
\subsection{Stochastic gradient numerical integrators}
For any numerical integrator for the underdamped Langevin dynamics \eqref{underdamped}, its stochastic version can be constructed by replacing the exact gradient $\nabla U(x)$ by $b(x,\omega_n)$, where the random variable $\omega_n$ is chosen independently at each iteration. 
In particular, the stochastic gradient versions of \eqref{EM} and \eqref{UBU} are constructed as follows.
\begin{definition}
Consider the following stochastic gradient numerical integrators.
\begin{enumerate}
\item \emph{Euler--Maruyama integrator.}
\begin{equation}
\left\{
\begin{aligned}
\tilde X_{n+1} & = \tilde X_n + h \tilde V_n, \\
\tilde V_{n+1} & = \tilde V_n - hb(\tilde X_n,\omega_n) - \gamma h \tilde V_n + \sqrt{2\gamma h} \xi_n,
\end{aligned}
\right.
\tag{SG-EM}
\label{SG-EM}
\end{equation}
where $\xi_n\sim \mathcal N(0,I_d)$ and $\omega_n\in\Omega$ are independently sampled at each time step. 
\item \emph{UBU integrator}. Consider two separate dynamics
\begin{equation*}
	\mathcal U:\left\{
	\begin{aligned}
	\dot x_t & = v_t, \\
	\dot v_t & = -\gamma v_t + \sqrt{2\gamma} \dot B_t,
	\end{aligned}
	\right.~~~~
	\tilde{\mathcal B}:\dot v_t = -b(x_t,\omega_n),
\end{equation*}
where both the $\mathcal U$ and $\tilde{\mathcal B}$ dynamics can be solved explicitly.
The update from $(\tilde X_n,\tilde V_n)$ to $(\tilde X_{n+1},\tilde V_{n+1})$ is given by
\begin{equation}
	(\tilde X_n,\tilde V_n)
	\xrightarrow{\mathcal U\mathrm{\,for\,}\frac h2}
	(\tilde Y_n,\tilde V_n^*)
	\xrightarrow{\tilde{\mathcal B}\mathrm{\,for\,} h}
	(\tilde Y_n,\tilde V_n^{**}) 
	\xrightarrow{\mathcal U\mathrm{\,for\,}\frac h2}
	(\tilde X_{n+1},\tilde V_{n+1}),
\tag{SG-UBU}
\label{SG-UBU}
\end{equation}
and $\omega_n\in\Omega$ is sampled independently at each time step.
\end{enumerate}
\end{definition}

The strong order of a stochastic gradient numerical integrator can be defined as follows.
\begin{definition}
\label{defi 2}
Suppose $p>0$ and $q\Ge 1$. The stochastic gradient numerical integrator 
\begin{equation*}
	\tilde Z_0 = (\tilde X_0,\tilde V_0) 
	\mapsto \tilde Z_1 = (\tilde X_1,\tilde V_1)
\end{equation*}
has strong order $(p,q)$ if its deterministic gradient version 
\begin{equation*}
	\tilde Z_0 = (\tilde X_0,\tilde V_0) \mapsto 
	\tilde Z_{0,1} = (\tilde X_{0,1}, \tilde V_{0,1})
\end{equation*}
has strong order $(p,q)$ in the sense of Definition~\ref{defi 1}, and
\begin{itemize}
\item The one-step update $\tilde Z_1$ is an unbiased approximation to $\tilde Z_{0,1}$ with respect to the random variable $\omega$, namely
\begin{equation}
	\mathbb E\big[
	\tilde Z_1 - \tilde Z_{0,1} \big| \tilde Z_{0,1}
	\big] = 0.
	\label{con 1}
\end{equation}
\item For some constant $C>0$, the difference $\tilde Z_1 - \tilde Z_{0,1}$ is bounded by
\begin{equation}
	\mathbb E|\tilde Z_1 - \tilde Z_{0,1}|^4 \Le Ch^4 (|\tilde Z_0|+1)^4.
	\label{con 2}
\end{equation}
\item For some constant $C>0$, the local error is bounded by
\begin{equation}
	\mathbb E|\tilde Z_1 - \tilde Z_0\{h\}|^4 \Le Ch^6 
	(|\tilde Z_0|+1)^4.
	\label{con 4}
\end{equation}
\item For any $r\Ge 1$, there exists a constant $C$ depending on $r$ such that
\begin{equation}
\mathbb E|\tilde Z_1|^r \Le C(|\tilde Z_0|+1)^r.
\label{con 5}
\end{equation}
\end{itemize}
A stochastic gradient numerical integrator has \emph{strong order $p$} if there exists $q\Ge 1$ such that this integrator has \emph{strong order $(p,q)$}.
\end{definition}
\begin{remark}
Remarks on the definition of the strong order:
\begin{enumerate}
\item The unbiased condition \eqref{con 1} characterizes the unbiased property of the numerical solution $\tilde Z_1$. The condition \eqref{con 1} is fulfilled by \eqref{SG-EM} and \eqref{SG-UBU}, but not necessarily by the stochastic gradient BAOAB, which involves the calculation of $b(x,\omega)$ twice in a single iteration.
\item In the condition \eqref{con 2}, the difference $\tilde Z_1 - \tilde Z_{0,1}$ is of order $\mathcal O(h)$ because the stochastic integral term vanishes, and thus $\tilde Z_1 - \tilde Z_{0,1}$ is determined by the difference between the full gradient and the stochastic gradient. Also, $\tilde Z_1 - \tilde Z_{0,1}$ is $\mathcal O(h)$ indicates the statistical error should be at least $\mathcal O(h^2\hspace{-1pt}+\hspace{-1pt}\frac1{Nh})$.
\item The conditions \eqref{con 4} and \eqref{con 5} can be viewed as the stochastic gradient versions of \eqref{cond 2} and \eqref{cond 3} in Definition~\ref{defi 1}. Utilizing the unbiased condition \eqref{con 1}, \eqref{con 5} directly implies
\begin{equation*}
	\mathbb E|\tilde Z_{0,1}|^r \Le C(|\tilde Z_0|+1)^r,
\end{equation*}
which is identical with \eqref{cond 3} in Definition~\ref{defi 1}. However, $\tilde Z_0(t)$ is not necessarily the expectation of $\tilde Z_0\{t\}$ with respect to  $\omega$, hence \eqref{con 4} does not imply \eqref{cond 2}.
\end{enumerate}
\end{remark}
Similar to the proof Lemma~\ref{lemma: strong order}, it is straightforward to verify the strong order of the stochastic gradient numerical integrators \eqref{SG-EM} and \eqref{SG-UBU}.
\begin{lemma}
\label{lemma: SG strong order}
Assume \textbf{A1}, \textbf{B1}. In the sense of Definition~\ref{defi 2}, 
\begin{enumerate}
\setzero
\item The stochastic gradient Euler--Maruyama integrator \eqref{SG-EM} has strong order $(1,1)$.
\item The stochastic gradient UBU integrator \eqref{SG-UBU} has strong order $(2,2)$.
\end{enumerate}
\end{lemma}
The key point of in the proof of Lemma~\ref{lemma: SG strong order} is to prove the unbiased condition \eqref{con 1} for \eqref{SG-EM} and \eqref{SG-UBU}. The proof is left in Appendix. 

Using the same strategies with Lemma~\ref{lemma: moment}, we can establish the uniform-in-time moment bound and the diffusion estimate for the stochastic gradient numerical integrators.
\begin{lemma}
\label{lemma: moment SG}
Assume \textbf{A1}, \textbf{A2}, \textbf{A4}. Suppose the stochastic gradient numerical integrator has strong order $(p,q)$, and let $r\Ge 1$ be an arbitrary constant. For some constant $\gamma_0>0$, whenever $\gamma\Ge \gamma_0$, there exist constants $C,h_0>0$ depending on $\gamma$ such that for any $0<h\Le h_0$,
\begin{equation}
	\sup_{n\Ge0} \mathbb E\big(
	|\tilde X_n| + |\tilde V_n| + 1
	\big)^{2r} \Le C \mathbb E\big(
		|\tilde X_0| + |\tilde V_0| + 1
		\big)^{2r}.
	\label{moment bound SG}
\end{equation}
\end{lemma}
\begin{lemma}
\label{lemma: diffusion SG}
Assume \textbf{A1}, \textbf{A2}, \textbf{A4}. Suppose the stochastic gradient numerical integrator has strong order $(p,q)$.
For some constant $\gamma_0>0$, whenever $\gamma\Ge \gamma_0$, there exist constants $C,h_0>0$ depending on $\gamma$ such that for any $0<h\Le h_0$,
\begin{equation}
	\mathbb E\big[|\tilde Z_N - \tilde Z_0|^2\big] \Le CNh \mathbb E
	(|\tilde Z_0|+1)^{2q}.
	\label{diffusion SG}
\end{equation}
\end{lemma}
The proof of Lemmas~\ref{lemma: moment SG} and \ref{lemma: diffusion SG} are provided in Appendix.
\subsection{Analysis of the statistical error}
To obtain the statistical error of the stochastic gradient numerical integrator, we need to evaluate the time average error $e(N,h)$, which is given by
\begin{equation*}
	e(N,h) = \frac1N \sum_{n=0}^{N-1}
	f(\tilde X_n,\tilde V_n) - \pi(f).
\end{equation*}
The randomness of $e(N,h)$ comes from the initial state $(\tilde X_0,\tilde V_0)$, the Brownian motion $(B_t)_{t\Ge0}$, and the random variables $(\omega_n)_{n=0}^{N-1}$.
The statistical error is given by $\mathbb E[e^2(N,h)]$.
\begin{theorem}
\label{theorem: SG}
Assume \textbf{A1}, \textbf{A2'}, \textbf{A3}, \textbf{A4}.
Suppose the stochastic gradient numerical integrator has strong order $(p,q)$. For some constant $\gamma_0>0$, whenever $\gamma\Ge \gamma_0$, there exist constants $C,h_0>0$ depending on $\gamma$ such that for any $0<h\Le h_0$,
\begin{equation*}
\mathbb E[e^2(N,h)] \Le C\bigg(
	h^{\min\{2p,2\}} + \frac1{Nh}
\bigg) \mathbb E\big(|\tilde X_0| + |\tilde V_0| + 1\big)^{4q}.
\end{equation*}
\end{theorem}
In particular, the statistical error of \eqref{SG-EM} and \eqref{SG-UBU} is both $\mathcal O\big(h^2\hspace{-1pt}+\hspace{-1pt}\frac1{Nh}\big)$.\\[6pt]
\emph{Proof of Theorem~\ref{theorem: SG}. Statistical error of the stochastic gradient numerical integrators.}\\[6pt]
Similar to the proof of Theorem~\ref{theorem: numerical error}, for integer each $n\Ge0$ we define the random variables
\begin{align*}
R_n & = \frac{\phi_h(\tilde Z_{n+1}) -\phi_h(\tilde Z_{n,1})}{h}, ~~~~ S_n = \frac{\phi_h(\tilde Z_{n,1}) - \phi_h(\tilde Z_n(h))}{h}, \\
T_n & = \frac{\phi_h(\tilde Z_n(h)) - \phi_h(\tilde Z_n)}{h}
+ f(\tilde Z_n) - \pi(f),
\end{align*}
where the definitions of $\tilde Z_n,\tilde Z_{n,1},\tilde Z_n(h)$ can be found in Section~\ref{section: 2}. Then, we have
\begin{equation}
	R_n + S_n + T_n = \frac{\phi_h(\tilde Z_{n+1}) - \phi_h(\tilde Z_n)}{h} + f(\tilde Z_n) - \pi(f).
	\label{RST}
\end{equation}
Summation of $R_n + S_n + T_n$ over $n$ in \eqref{RST} yields a similar expression with \eqref{rountine 2}:
\begin{equation}
e(N,T) = \frac{\phi_h(\tilde Z_0) - \phi_h(\tilde Z_N)}{Nh} + \frac1N \sum_{n=0}^{N-1} (R_n + S_n + T_n),
\label{e SG}
\end{equation}
and the only difference is the additional term $R_n$ due to stochastic gradients.
Similar to the proof of Theorem~\ref{theorem: numerical error convex}, the quantities $S_n$ and $T_n$ satisfy
\begin{align}
\mathbb E\Bigg[
\bigg|\sum_{n=0}^{N-1} S_n \bigg|^2
\Bigg] & \Le CN^2 \bigg(h^{2p} + \frac1{Nh}\bigg) \mathbb E(|\tilde Z_0|+1)^{4q}, \label{final S} \\
\mathbb E\Bigg[
\bigg|\sum_{n=0}^{N-1} T_n \bigg|^2
\Bigg] & = \sum_{n=0}^{N-1} \mathbb E[T_n^2] \Le 
\frac{CN}h \mathbb E(|\tilde Z_0|+1)^2.
\label{final T}
\end{align}
Here we recall that the numerical integrator has strong order $(p,q)$, and the random variables $\{T_n\}_{n=0}^{N-1}$ are unbiased and mutually independent.

Now, we focus on estimating $R_n$. Notably, $R_n$ can be expressed as
\begin{equation*}
	R_n = \frac1h\big(\tilde Z_{n+1} - \tilde Z_{n,1}\big)
	\cdot \int_0^1 \nabla \phi_h\big(\theta
	\tilde Z_{n+1} + (1-\theta) \tilde Z_{n,1}
	\big) \d\theta,
\end{equation*}
which can be decomposed into $R_{n,1} + R_{n,2}$, where
\begin{align*}
R_{n,1} & = \frac1h (\tilde Z_{n+1} - \tilde Z_{n,1})\cdot 
\nabla \phi_h(\tilde Z_{n,1}), \\
R_{n,2} & = \frac1h (\tilde Z_{n+1} - \tilde Z_{n,1})\cdot 
\int_0^1 \Big(\nabla \phi_h \big(\theta \tilde Z_{n+1} + (1-\theta) \tilde Z_{n,1}\big) 
- \nabla \phi_h(\tilde Z_{n,1})\Big)
\d\theta.
\end{align*}
From the condition \eqref{cond 1}, the random variables $\{R_{n,1}\}_{n=0}^{N-1}$ are unbiased and thus mutually independent. Utilizing the boundedness of $\nabla \phi_h(x,v)$ and $\nabla^2 \phi_h(x,v)$ as stated in Lemmas~\ref{lemma: phi h 1} and \ref{lemma: phi h 2}, we derive
\begin{equation*}
|R_{n,1}| \Le \frac Ch |\tilde Z_{n+1} - \tilde Z_{n,1}|, ~~~~
|R_{n,2}| \Le \frac Ch |\tilde Z_{n+1} - \tilde Z_{n,1}|^2.
\end{equation*}
Thus, by the condition \eqref{con 2} and Lemma~\ref{lemma: moment SG}, we obtain
\begin{align}
\mathbb E\Bigg[
\bigg|\sum_{n=0}^{N-1} R_{n,1}\bigg|^2
\Bigg] & = \sum_{n=0}^{N-1} \mathbb E\big[R_{n,1}^2\big]
\Le \frac{C}{h^2} \sum_{n=0}^{N-1} \mathbb E\big[|\tilde Z_{n+1} - \tilde Z_{n,1}|^2\big] \notag \\
& \Le C \sum_{n=0}^{N-1}\mathbb E(|\tilde Z_n|+1)^2 \Le CN
\mathbb E(|\tilde Z_0|+1)^2.
\label{R estimate 1}
\end{align}
By applying Cauchy's inequality and the condition \eqref{con 2}, we further derive
\begin{align}
\mathbb E\Bigg[
\bigg|\sum_{n=0}^{N-1} R_{n,2}\bigg|^2
\Bigg] & \Le N \sum_{n=0}^{N-1} \mathbb E\big[R_{n,2}^2\big]  \Le 
\frac{CN}{h^2} \sum_{n=0}^{N-1} \mathbb E\big[|\tilde Z_{n+1} - \tilde Z_{n,1}|^4\big] \notag \\
& \Le \frac{CN}{h^2} \sum_{n=0}^{N-1} 
 h^4 \mathbb E(|\tilde Z_n|+1)^4 \Le CN^2 h^2 \mathbb E(|\tilde Z_0|+1)^4.
\label{R estimate 2}
\end{align}
Combining \eqref{R estimate 1} and \eqref{R estimate 2}, we conclude
\begin{equation}
\mathbb E\Bigg[
\bigg|\sum_{n=0}^{N-1} R_{n}\bigg|^2
\Bigg] \Le CN^2\bigg(h^2 + \frac1{Nh}\bigg) \mathbb E(|\tilde Z_0|+1)^4.
\label{final R}
\end{equation}

Finally, from the expression of $e(N,h)$ in \eqref{e SG}, the diffusion estimate \eqref{diffusion SG} in Lemma~\ref{lemma: diffusion SG}, and the quantitative estimates of $R_n,S_n,T_n$ in \eqref{final S}\eqref{final T}\eqref{final R}, we obtain
\begin{align*}
 \mathbb E[e^2(N,h)] & \Le 
\frac{2}{N^2h^2}\mathbb E\big[|\phi_h(\tilde Z_0) - \phi_h(\tilde Z_N)|^2\big] + \frac2{N^2} \mathbb E\Bigg[
\bigg|\sum_{n=0}^{N-1} (R_n + S_n + T_n ) \bigg|^2
\Bigg] \\
& \Le \frac{C}{N^2h^2} \mathbb E\big[|\tilde Z_0 - \tilde Z_N|^2\big] + \frac6{N^2}
\mathbb E\Bigg[
\bigg|\sum_{n=0}^{N-1} R_{n}\bigg|^2 + 
\bigg|\sum_{n=0}^{N-1} S_{n}\bigg|^2 + 
\bigg|\sum_{n=0}^{N-1} T_{n}\bigg|^2
\Bigg]
\\
& \Le C\bigg(h^{2p} + \frac1{Nh}\bigg) \mathbb E(|\tilde Z_0|+1)^{4q} + 
C\bigg(h^2 + \frac1{Nh}\bigg) \mathbb E(|\tilde Z_0|+1)^{4q} \\
& \Le C\bigg(h^{\min\{2p,2\}} + \frac1{Nh}\bigg) \mathbb E(|\tilde Z_0|+1)^{4q},
\end{align*}
which completes the proof of Theorem~\ref{theorem: SG}. \hfill $\blacksquare$
\begin{remark}
In the proof of Theorem~4, the error term $\mathcal O(h^2)$ in the statistical error comes from $\mathbb E\big[|R_{n,2}|^2\big]$, which is bounded by $\mathbb E\big[|\tilde Z_{n+1} - \tilde Z_{n,1}|^4\big]$. If we make use the mini-batch form of the stochastic gradient $b(x,\omega)$ in \eqref{mini batch}, the upper bound of the statistical error $\mathbb E[e^2(N,h)]$ should be explicit on the batch size $\mathcal B$.
\end{remark}
\section{Numerical verification of the statistical error}
\label{section: 7}
To illustrate the statistical error of the numerical integrators presented in this paper, we employ a straightforward example focusing on the dependence on the time step $h$. Let's consider the 1-dimensional potential function
\begin{equation*}
	U(x) = \frac12|x|^2 + \sin x,~~~~x\in\mathbb R^1,
\end{equation*}
with gradient $\nabla U(x) = x + \cos(x)$. We construct the stochastic gradient as
\begin{equation*}
	b(x,\omega) = \omega_1x + \omega_2 + \cos x,
\end{equation*}
where $\omega = (\omega_1,\omega_2)$ is a random variable in $\mathbb{R}^2$ such that
$$
\omega_1\sim\mathrm{Uniform}(0.2,1.8),~~~~
\omega_2\sim \mathcal N(0,0.4^2).
$$
Clearly, $b(x,\omega)$ serves as an unbiased approximation to $\nabla U(x)$.

In the numerical simulations, we set the damping rate $\gamma = 2$, the initial state $(X_0,V_0) = (0.2,-0.3)$, the total simulation time $T = 10^7$, and the time step $h$ varies from $2^{-1}$ to $2^{-6}$. We employ 100 independent trajectories to compute the expectation in $\mathbb{E}[e^2(N,h)]$. The statistical average $\mathbb{E}[e^2(N,h)]$ for the numerical integrators \eqref{EM}, \eqref{UBU} and their stochastic gradient versions \eqref{SG-EM}, \eqref{SG-UBU} is depicted in Figure 1.
\begin{center}
\includegraphics[width=240pt]{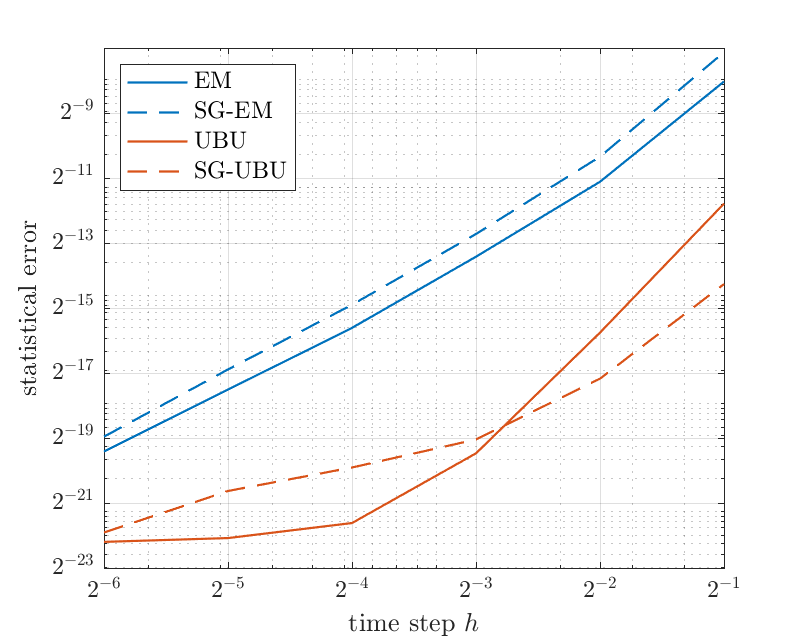}
\captionof{figure}{The statistical error $\mathbb E[e^2(N,h)]$ for the numerical integrators \eqref{EM}, \eqref{UBU} and their stochastic gradient versions \eqref{SG-EM}, \eqref{SG-UBU} for various $h$.} 
\end{center}
Observations reveal that \eqref{UBU} and \eqref{SG-UBU} exhibit smaller statistical variances compared to the first-order integrator \eqref{EM} and \eqref{SG-EM}. Particularly, \eqref{UBU} demonstrates a swift decline in the statistical error as the time step $h$ diminishes. The reason that the statistical error reaches a plateau as the time step $h$ decreases is that to the dominance of the $\mathcal{O}\big(\frac1{Nh}\big)$ component in the error term.
\section{Conclusion}
\label{section: 8}
In this paper we employ the discrete Poisson solution $\phi_h(x,v)$ to study the statistical error of the numerical integrators for the underdamped Langevin dynamics \eqref{underdamped}. When the potential function $U(x)$ is strongly convex outside a ball, we prove that the statistical error is $
\mathcal O(h^{2p-1}\hspace{-1pt}+\hspace{-1pt}\frac1{Nh})$, where $p>0$ is the strong order of the numerical integrator. When $U(x)$ is strongly convex in $\mathbb R^d$, the estimate is enhanced to be $\mathcal O(h^{2p}\hspace{-1pt}+\hspace{-1pt}\frac1{Nh})$. For the stochastic gradient numerical integrator, when $U(x)$ is strongly convex in $\mathbb R^d$, we show that the statistical error is $\mathcal O(h^2 \hspace{-1pt} + \hspace{-1pt} \frac1{Nh})$. Our results provide the estimates of the statistical error for a broad class of numerical integrators for the underdamped Langevin dynamics \eqref{underdamped}, and relaxes the constraint on the time step $h$.

There are several potential avenues for future research.
Firstly, for stochastic gradient numerical integrators, it would be intriguing to derive an explicit bound of the statistical error with clear dependence on the batch size $\mathcal B$ and the spatial dimension $d$.
Secondly, investigating the stochastic variance reduced gradient (SVRG) dynamics \cite{SGLD_1,UBU_3} could be worthwhile.

\section*{Acknowledgments}
The work of Z. Zhou is partially supported by the National Key R\&D Program of China
(Project No. 2020YFA0712000, 2021YFA1001200), and the National Natural Science Foundation
of China (Grant No. 12031013, 12171013).

The authors would like to thank Jian-Guo Liu (Duke), Benedict Leimkuhler (UoE), Daniel Paulin (UoE) and Peter Whalley (UoE) for the helpful discussions.

\printbibliography

\section*{Appendix~~Addtional proofs in Section~\ref{section: 3} to Section~\ref{section: 6}}
\emph{Proof of Lemma \ref{lemma: strong order}. Strong order estimate of various numerical integrators.}\\[6pt]
The proof for the Euler--Maruayma integrator is standard and can be found in \cite{Lyapunov}.
The proof for the UBU integrator \eqref{UBU} is the same with Section 7.6 of \cite{UBU_2}, and the only difference is that the error terms require the fourth moments rather than second.
Being consistent with the notations in \cite{UBU_2}, introduce
\begin{equation*}
\mathcal E(t) = \exp(-\gamma t) \Le 1,~~~~
\mathcal F(t) = \int_0^t \mathcal E(s)\d s = \frac{1-e^{-\gamma t}}{\gamma} \Le t.
\end{equation*}
Then the UBU integrator from $(X_0,V_0)$ to $( X_1, V_1)$ can be equivalent written as
\begin{subequations}
\begin{align}
 V_1 & = \mathcal E(h) V_0 - h \mathcal E(h/2)
\nabla U(Y_0) + \sqrt{2\gamma} \int_0^h
\mathcal E(h-s)\d B_s, \\
 X_1 & = X_0 + \mathcal F(h) V_0 - h \mathcal F(h/2)
\nabla U(Y_0) + \sqrt{2\gamma} \int_0^h
\mathcal F(h-s)\d B_s, \\
Y_0 & = X_0 + \mathcal F(h/2) V_0 + \sqrt{2\gamma}
\int_0^{h/2} \mathcal F(h/2-s)\d B_s.
\end{align}
\label{UBU explicit}
\end{subequations}
Denote the intermediate difference
\begin{equation}
	\Delta_y = X_0(h/2) - Y_0 = -\int_0^{h/2} 
	\mathcal F(h/2-s)\nabla U(X_0(s))\d s.
\end{equation}
then by \textbf{A1}, the difference $\Delta_y$ has the estimate
\begin{equation*}
	\mathbb E|\Delta_y|^4 \Le Ch^4 \mathbb E\Bigg[\bigg(\int_0^{h/2} 
	(|X_0(s)|+1)\d s\bigg)^4\Bigg] \Le Ch^8
	(|Z_0|+1)^4.
\end{equation*}
According to Equation (40) of \cite{UBU_2}, $\Delta_v$ can be written as
\begin{equation}
\Delta_v = h\mathcal E(h/2) \big(
\nabla U(X_0(h/2)) - \nabla U(Y_0)
\big) + I_1 + I_2 + I_3 + I_4 + I_5.
\label{Delta v}
\end{equation}
The quantities $I_1,I_2,I_3,I_4$ are the high-order terms, and satisfy
\begin{equation*}
\mathbb E|I_1|^4,~\mathbb E|I_2|^4,~
\mathbb E|I_3|^4,~\mathbb E|I_4|^4 \Le C h^{12}(|Z_0|+1)^{8}.
\end{equation*}
The quantity $I_5$ is a stochastic integral given by
\begin{equation*}
I_5  = \sqrt{2\gamma}
\int_{h/2}^h\d s \int_{h/2}^s\d s' \int_{h-s'}^{s'}
\mathcal E(h-s'') \nabla^2 U(X_0(s''))\d B_{s''},
\end{equation*}
which simply satisfies $\mathbb E|I_5|^4 \Le Ch^{10}$. Therefore, $I_5$ is the mean-zero random variable in the local error $\Delta_v$, and the remaining terms in \eqref{Delta v} are the high-order term.

According to Equation (41) of \cite{UBU_2}, $\Delta_x$ can be written as
\begin{equation}
\Delta_x = h\mathcal F(h/2) \big(
\nabla U(X_0(h/2)) - \nabla U(Y_0)
\big) + I_6 + I_7,
\label{Delta_x}
\end{equation}
where the quantities $I_6,I_7$ satisfy
$$
\mathbb E|I_6|^4,\mathbb E|I_7|^4 \Le Ch^{12}(|Z_0|+1)^8.
$$
Hence from \eqref{Delta_x} $\mathbb E|\Delta_x|^4$ is bounded by
\begin{equation*}
\mathbb E|\Delta_x|^4 \Le Ch^4 \mathbb E|\Delta_y|^4 + Ch^{12}(|Z_0|+1)^8 \Le 
Ch^{12}(|Z_0|+1)^8,
\end{equation*}
In conclusion, the UBU integrator \eqref{UBU} has strong order $(2,2)$.

Finally, we verify the linear stability condition \eqref{cond 2} for \eqref{UBU}. Note that
\begin{subequations}
\begin{align}
&  V_1 - V_0 - \frac{h}2 \nabla U(X_0) - \sqrt{2\gamma} B_h = 
h\mathcal E(h/2)\big(\nabla U(X_0)-\nabla U(Y_0)\big)\notag  \\
& \hspace{3cm} + \sqrt{2\gamma} \int_0^h \big(\mathcal E(h-s)-1\big)\d B_s +
\mathcal O(h^2(|Z_0|+1)), \\
& X_1 - X_0 - hV_0  = -h\mathcal F(h/2) \nabla U(Y_0) + \sqrt{2\gamma}\int_0^h \mathcal F(h-s)\d B_s + \mathcal O(h^2(|Z_0|+1)).
\end{align}
\label{X1V1 expression}
\end{subequations}
Using the estimates
\begin{align*}
\mathbb E\big|\nabla U(Y_0) - \nabla U(X_0)\big|^4 & \Le C
\mathbb E|Y_0 - X_0|^4 \Le Ch^4(|Z_0|+1)^4, \\
\mathbb E|\nabla U(Y_0)|^4 & \Le C \mathbb E(|Y_0|+1)^4 \Le 
C (|Z_0|+1)^4
\end{align*}
in \eqref{X1V1 expression}, we immediately obtain
\begin{align*}
\mathbb E\Big|V_1 - V_0 - \frac{h}2 \nabla U(X_0) - \sqrt{2\gamma} B_h\Big|^4 & \Le Ch^6(|Z_0|+1)^4, \\
\mathbb E\big|X_1 - X_0 - hV_0\big|^4 & \Le Ch^6(|Z_0|+1)^4.
\end{align*}
As a consequence, we have $\mathbb E|Z_1 - Z_0(h)|^4 \Le Ch^6(|Z_0|+1)^4$.
\hfill $\blacksquare$\\[6pt]
\emph{Proof of Lemma~\ref{lemma: moment}. Uniform-in-time moment estimates of the numerical solution.}\\[6pt]
We employ the same procedure with Theorem 7.2 of \cite{Lyapunov} to prove this result. The generator of the underdamped Langevin dynamics \eqref{underdamped} is clearly
\begin{equation*}
\mathcal L = v\cdot \nabla_x - (\nabla U(x) + \gamma v)
\cdot \nabla_v + \gamma \Delta_v.
\end{equation*}
To proceed, we prove that $\mathcal H(x,v)$ defined in \eqref{Lyapunov} is a Lyapunov function of \eqref{underdamped},
namely, there exist constants $a,b>0$ such that
\begin{equation}
	(\mathcal L \mathcal H)(x,v) \Le -a\mathcal H(x,v) + b,
	~~~~\forall x,v \in \mathbb R^d.
	\label{Lyapunov condition 0}
\end{equation}
Since $x\cdot \nabla U(x) \Ge m|x|^2 - C_0$ according to \eqref{U diffu}, we obtain
\begin{align*}
(\mathcal L\mathcal H)(x,v) & = v\cdot\nabla_x \mathcal H - (\nabla U(x) + \gamma v)\cdot \nabla_v \mathcal H + \gamma \Delta_v \mathcal H \\
& = -2\Big(
(\gamma-1) |v|^2 + x^\T \nabla U(x) + v^\T \nabla U(x)
\Big) + 2\gamma \\
& \Le -2\Big(
(\gamma-1)|v|^2 + m|x|^2 - C_1|v|(|x|+1) - C_0
\Big) + 2\gamma \\
& \Le -2\Big((\gamma-2)|v|^2 + m|x|^2 - C_1|x||v|\Big) + C.
\end{align*}
By choosing $\gamma$ sufficiently large, there exists a constant $a>0$ depending on $\gamma$ such that
\begin{equation*}
	2\Big((\gamma-2)|v|^2 + m|x|^2 - C_1|x||v|\Big) \Ge 
	a\Big(
	\gamma|x|^2 + 2x^\T v + |v|^2
	\Big),~~~~
	\forall x,v\in\mathbb R^d,
\end{equation*}
and thus \eqref{Lyapunov condition 0} holds true.
For notational convenience, denote
\begin{equation*}
	z = (x,v),~~~~Z_0 = (X_0,V_0),~~~~
	Z_0(t) = (X_0(t),V_0(t)),~~~~
	Z_1 = (X_1,V_1).
\end{equation*}
Since $r\Ge1$, $\mathcal H_r(z) := \big(\mathcal H(z)+1\big)^r$ is a Lyapunov function of \eqref{underdamped} according to Lemma 3.3 of \cite{Lyapunov}, and thus for some constants $a,b>0$, the Lyapunov condition holds:
\begin{equation}
	(\mathcal L\mathcal H_r)(z) \Le -a\mathcal H_r(z) + b,
	~~~~\forall z = (x,v) \in \mathbb R^{2d}.
	\label{Lyapunov condition}
\end{equation}
Now we prove for some constants $a_1,b_1>0$, the numerical Lyapunov condition holds:
\begin{equation}
	\mathbb E\big[\mathcal H_r(Z_1)\big] \Le 
	e^{-a_1 h} \mathcal H_r(Z_0) + b_1h,
	\label{Lyapunov condition numerical}
\end{equation}
where $Z_0 = (X_0,V_0)$ is the fixed initial state, and $\mathcal Z_1 = (X_1,V_1)$ is the one-step update of the numerical integrator. Using the Lyapunov condition \eqref{Lyapunov condition}, we have
\begin{align}
\mathbb E\big[\mathcal H_r(Z_1)\big] & \Le 
\mathbb E\big[\mathcal H_r(Z_0(h))\big] + 
\mathbb E\big|\mathcal H_r(Z_1) - \mathcal H_r(Z_0(h))\big| \notag \\
& \Le e^{-ah} \mathcal H_r(Z_0) + 
\frac{b}a(1-e^{-ah}) + \mathbb E\int_0^1
|Z_1 - Z_0(h)| \big|
\nabla\mathcal H_r\big(sZ_1 + (1-s)Z_0(h)\big)
\big|\d s \notag \\
& \Le e^{-ah} \mathcal H_r(Z_0) + 
\frac{b}a(1-e^{-ah}) + C\,\mathbb E\Big[\big|Z_1 - Z_0(h)\big| 
\big(|Z_1| + |Z_0(h)| + 1\big)^{2r-1} \Big] \notag \\
& \Le e^{-ah} \mathcal H_r(Z_0) + 
\frac{b}a(1-e^{-ah}) + C\sqrt{\mathbb E\big|Z_1 - Z_0(h)\big|^2}
\sqrt{\mathbb E\big(|Z_1| + |Z_0(h)| + 1\big)^{4r-2}}\notag  \\
& \Le e^{-ah} \mathcal H_r(Z_0) + 
\frac{b}a(1-e^{-ah}) + C (|Z_0|+1)^{2r} h^{\frac32} \notag \\
& \Le \big(e^{-ah} + Ch^{\frac32}\big)\mathcal H_r(Z_0) + 
\frac{b}a(1-e^{-ah}).
\label{Lyapunov derivation}
\end{align}
By choosing $h$ sufficiently small we have
\begin{equation*}
	e^{-ah} + Ch^{\frac32} \Le  e^{-\frac12ah},
\end{equation*}
and thus \eqref{Lyapunov derivation} yields the inequality
\begin{equation}
\mathbb E\big[\mathcal H_r(Z_1)\big] \Le 
e^{-\frac12 ah} \mathcal H_r(Z_0) + \frac ba(1-e^{-ah}).
\label{12 inequality}
\end{equation}
which is exactly the desired result \eqref{Lyapunov condition numerical}.

For the numerical solution $Z_n = (X_n,V_n)$ evolved by the numerical integrator, the numerical Lyapunov condition
\eqref{Lyapunov condition numerical} ensures
\begin{equation*}
	\sup_{n\Ge 0} \mathbb E\big[\mathcal H_r(Z_n)\big]
	\Le \max\bigg\{
	\frac{2b_1}{a_1}, \mathcal H_r(Z_0)
	\bigg\} \Le C \mathbb E\big(|Z_0| + 1\big)^{2r},
\end{equation*}
producing the uniform-in-time moment bound \eqref{uniform moment bound}.\\[6pt]
\emph{Proof of Lemma~\ref{lemma: diffusion}. Diffusion estimate of the numerical integrator.}\\[6pt]
We prove that for given initial state $Z_0 = (X_0,V_0)$, the error $Z_1 - Z_0$ can be written as
\begin{equation}
	Z_1 - Z_0 = \bar M_h + \bar N_h,
	\label{sum_1}
\end{equation}
where $\bar M_h$ is a mean-zero random variable, $\bar N_h$ is a high-order term, and
\begin{equation}
	\mathbb E\big[|\bar M_h|^2\big] \Le Ch (|Z_0|+1)^{2q},
	~~~~
	\mathbb E\big[|\bar N_h|^2\big] \Le Ch^2 (|Z_0|+1)^{2q}.
	\label{sum_2}
\end{equation}
On the one hand, according to \eqref{cond 1}, the local error $Z_1 - Z_0(h)$ can be written as $M_h + N_h$. On the other hand, $Z_0(h) - Z_0$ is the summation of a mean-zero random variable and a high-order term from It\^o's calculus, hence \eqref{sum_1} and \eqref{sum_2} hold true.

Now for each $n\Ge0$, $Z_{n+1} - Z_n$ can be written as $\bar M_{n,h}$ and $\bar N_{n,h}$, where $\bar M_{n,h}|Z_n$ is a mean-zero random variable, and $\bar N_{n,h}|Z_n$ is a high-order term. Also, from \eqref{sum_2} we obtain
\begin{equation*}
\mathbb E\big[
|\bar M_{n,h}|^2
\big] \Le Ch \mathbb E(|Z_n|+1)^{2q},~~~~
\mathbb E\big[
|\bar N_{n,h}|^2 
\big] \Le Ch^2 \mathbb E(|Z_n|+1)^{2q}.
\end{equation*}
Using the uniform-in-time $2q$-th moment estimate in Lemma~\ref{lemma: moment}, we have
\begin{equation*}
\mathbb E\big[
|\bar M_{n,h}|^2
\big] \Le Ch \mathbb E(|Z_0|+1)^{2q},~~~~
\mathbb E\big[
|\bar N_{n,h}|^2 
\big] \Le Ch^2 \mathbb E(|Z_0|+1)^{2q}.
\end{equation*}
Since the random variables $\{\bar M_{n,h}\}_{n=0}^{N-1}$ are mutually independent, we obtain
\begin{align*}
\mathbb E\big[|Z_N - Z_0|^2\big] & = \mathbb E\Bigg[
\bigg|\sum_{n=0}^{N-1} (\bar M_{n,h} + \bar N_{n,h})\bigg|^2\Bigg] \\
& \Le 2\mathbb E\Bigg[\bigg|\sum_{n=0}^{N-1}
\bar M_{n,h}
\bigg|^2\Bigg] + 2\mathbb E\Bigg[\bigg|\sum_{n=0}^{N-1}
\bar N_{n,h}
\bigg|^2\Bigg] \\
& \Le 2\sum_{n=0}^{N-1} \mathbb E\big[|\bar M_{n,h}|^2\big] + 2 N \sum_{n=0}^{N-1} \mathbb E\big[|\bar N_{n,h}|^2\big] \\
& \Le C(Nh + N^2h^2) \mathbb E(|Z_0|+1)^{2q}.
\end{align*}
Hence if $Nh\Le1$, we obtain the desired result $\mathbb E\big[|Z_N - Z_0|^2\big] \Le CNh \mathbb E(|Z_0|+1)^{2q}$.

Otherwise when $Nh>1$, utilizing Lemma~\ref{lemma: moment} we derive
\begin{equation*}
	\mathbb E\big[|Z_N - Z_0|^2\big] 
	\Le 2 \mathbb E\big[|Z_N|^2 + |Z_0|^2\big] \Le C Nh \mathbb E(|Z_0|+1)^2,
\end{equation*}
which completes the proof of Lemma~\ref{lemma: diffusion}.
\hfill $\blacksquare$ \\[6pt]
\emph{Proof of Lemma~\ref{lemma: estimate u 1}. Estimate $u(x,v,t)$ under strong convexity outside a ball.}\\[6pt]
For any $(x,v)\in\mathbb R^{2d}$, we have
\begin{equation*}
\big|(e^{t\mathcal L}f)(x,v) - \pi(f)\big| = 
\big|(\delta^{x,v}p_t)(f) - \pi(f)\big| \Le C \mathcal W_1(\delta^{x,v} p_t,\pi).
\end{equation*}
Using the Wasserstein-1 contractivity in Theorem~\ref{theorem: contractivity}, we have
\begin{equation*}
\mathcal W_1(\delta^{x,v}p_t,\pi) \Le Ce^{-\lambda t}
\mathcal W_1(\delta^{x,v},\pi) \Le C e^{-\lambda t} (|x|+|v|+1).
\end{equation*}
Hence from the two inequalities above we obtain the bound of $|u(x,v,t)|$,
\begin{equation*}
|u(x,v,t)| =
\big|(e^{t\mathcal L}f)(x,v) - \pi(f)\big| \Le 
C e^{-\lambda t} (|x|+|v|+1).
\end{equation*}

To estimate $\nabla u(x,v,t)$, pick two initial states $(x,v)$ and $(x',v')$ in $\mathbb R^{2d}$. Again using the Wasserstein-1 contractivity, we have
\begin{align*}
|u(x,v,t) - u(x',v',t)| 
	& = \big|(e^{t\mathcal L}f)(x,v) -
	(e^{t\mathcal L}f)(x',v')\big| = 
	\big|(\delta^{x,v} p_t)(f) - (\delta^{x',v'} p_t)(f)\big| \\
	& \Le \mathcal W_1\big(\delta^{x,v}p_t,\delta^{x',v'}p_t\big) \Le 
	C  e^{-\lambda t} \mathcal W_1\big(
	\delta^{x,v}, \delta^{x',v'}
	\big) \\
	& = C e^{-\lambda t}
	(|x-x'|+|v-v'|).
\end{align*}
As $|x-x'|+|v-v'|\rightarrow0$, we obtain the desired result.  \hfill $\blacksquare$\\[6pt]
\emph{Proof of Lemma~\ref{lemma: qp decay}. Exponential decay of the tangent processes $(Q_t,P_t)$.} Denote
\begin{equation*}
	W_t = \begin{pmatrix}
	Q_t \\ P_t
	\end{pmatrix} \in \mathbb R^{2d\times d},~~~
	S = \begin{pmatrix}
	\gamma I_d & I_d \\
	I_d & I_d
	\end{pmatrix} \in \mathbb R^{2d\times 2d},~~~
	A_t = \begin{pmatrix}
	O_d & -I_d \\
	\nabla^2 U(x_t) & \gamma I_d
	\end{pmatrix} \in \mathbb R^{2d\times 2d}.
\end{equation*}
Then the Lyapunov function $\mathcal H(Q_t,P_t)$ can be shortly written as $\mathcal H(W_t) = \Tr[W_t^\T S W_t]$.
From the strong convexity condition in \textbf{A2} we have 
\begin{equation*}
	m I_d \Prec \nabla^2 U(x_t) \Prec C_1 I_d,~~~~
	\mathrm{for~any~}t\Ge0.
\end{equation*}
Note that for any $t\Ge0$, there holds the inequality
\begin{align*}
	SA_t + A_t^\T S & = \begin{pmatrix}
	2\nabla^2 U(x_t) & \nabla^2 U(x_t) \\
	\nabla^2 U(x_t) & 2(\gamma-1) I_d
	\end{pmatrix} \\
	& \Succ
	\begin{pmatrix}
	\nabla^2 U(x_t) & O_d \\
	O_d & 2(\gamma-1) I_d - \nabla^2 U(x_t)
	\end{pmatrix} \Succ 
	\begin{pmatrix}
	m I_d & O_d \\
	O_d & \big(2(\gamma-1) - C_1\big) I_d
	\end{pmatrix}.
\end{align*}
When the damping rate $\gamma\Ge C_1+1$, there exists a constant $\lambda>0$ such that
\begin{equation}
SA_t + A_t^\T S \Succ
\begin{pmatrix}
m I_d & O_d \\
O_d & C_1 I_d
\end{pmatrix} \Succ \lambda S,
\label{SA bound}
\end{equation}
and thus the Lyapunov function $\mathcal H(W_t)$ satisfies
\begin{equation*}
\frac{\d}{\d t}\mathcal H(W_t) = 
\frac{\d}{\d t}\Tr[W_t^\T S W_t] = \Tr[W_t^\T (SA_t + A_t^\T S) W_t] \Le -\lambda \Tr[ W_t^\T S W_t] = -\lambda \mathcal H(W_t).
\end{equation*}
As a consequence, $\mathcal H(W_t)$ has exponential decay, and
\begin{equation*}
\mathcal H(W_t) \Le e^{-\lambda t} \mathcal H(W_0) \Le 
C e^{-\lambda t} (|Q_0|+|P_0|)^2,
\end{equation*}
which completes the proof.
 \hfill $\blacksquare$\\[6pt]
\emph{Proof of Lemma~\ref{lemma: qp diff decay}. Coupling error the tangent processes $(Q_t,P_t)$ and $(Q_t',P_t')$.}\\[6pt]
Similar to the proof of Lemma~\ref{lemma: qp decay}, denote
\begin{equation*}
W_t = \begin{pmatrix}
Q_t \\ P_t
\end{pmatrix},~~~~
W_t' = \begin{pmatrix}
Q_t' \\ P_t'
\end{pmatrix},~~~~
A_t = \begin{pmatrix}
O_d & -I_d \\ \nabla^2 U(x_t) & \gamma I_d
\end{pmatrix},~~~~
A_t' = \begin{pmatrix}
O_d & -I_d \\ \nabla^2 U(x_t') & \gamma I_d
\end{pmatrix},
\end{equation*} 
where the matrices $\nabla^2 U(x_t)$ and $\nabla^2 U(x_t')$ satisfy
\begin{equation*}
m I_d \Prec \nabla^2 U(x_t) , \nabla^2 U(x_t') \Prec C_1 I_d,~~~~
\mathrm{for~any~}t\Ge0.
\end{equation*}
Let $\delta = |x - x'| + |v-v'|$ be the deviation in the initial states.

The strong solutions $(x_t,v_t)$ and $(x_t',v_t')$ are driven by the same Brownian motion (synchronous coupling), hence from the proof of Theorem~1 of \cite{coupling_1}, there exists a constant $\gamma_0\Ge C_1+1$ such that the difference between $(x_t,v_t)$ and $(x_t',v_t')$ is bounded by
\begin{equation*}
	|x_t - x_t'| + |v_t - v_t'| \Le C\delta,~~~~
	\mathrm{for~any~}t\Ge0.
\end{equation*}
Write $B_t = \nabla^2 U(x_t) - \nabla^2 U(x_t') \in \mathbb R^{d\times d}$. Since $\nabla^3 U(x)$ is bounded by \textbf{A1}, we have
\begin{equation*}
	|B_t| \Le |\nabla^2 U(x_t) - \nabla^2 U(x_t')|
	\Le C|x_t-x_t'| \Le C\delta,~~~~
	\forall t\Ge0.
\end{equation*}
The time derivative of the Lyapunov function $\mathcal H_t := \mathcal H(W_t - W_t')$ is given by
\begin{align*}
\frac{\d}{\d t}\mathcal H_t & = \frac{\d}{\d t}
\Tr[(W_t-W_t')^\T S (W_t - W_t')] \notag \\
& = - 2 \Tr[(W_t - W_t') S (A_t W_t - A_t' W_t')] \notag \\
& = - 2 \Tr[(W_t - W_t')^\T SA_t' (W_t - W_t')] - 
2 \Tr[(W_t - W_t')^\T S (A_t - A_t') W_t].
\end{align*}
According to \eqref{SA bound}, for some constant $\lambda_1>0$ there is $S A_t' + (A_t')^\T S \Succ \lambda_1 S$, 
and thus 
\begin{equation}
\frac{\d}{\d t}\mathcal H_t \Le -\lambda_1 
\mathcal H_t - 2\Tr[(W_t - W_t')^\T S B_t W_t].
\label{F_t derivative}
\end{equation}
From Lemma~\ref{lemma: qp decay}, there exists a constant $\lambda_2>0$ such that $|W_t| \Le C e^{-\lambda_2 t}(|Q_0|+|P_0|)$, then
\begin{equation}
\big|(W_t-W_t')^\T S B_t W_t\big| \Le 
C|B_t||W_t-W_t'||W_t| \Le 
C e^{-\lambda_2 t}
\sqrt{\mathcal H_t}
(|Q_0|+|P_0|)\delta.
\label{ww estimate}
\end{equation}
We derive from \eqref{F_t derivative} and \eqref{ww estimate} that
\begin{equation}
\frac{\d}{\d t} \mathcal H_t \Le 
-\lambda_1 \mathcal H_t + C e^{-\lambda_2t}
\sqrt{\mathcal H_t} (|Q_0|+|P_0|)\delta.
\end{equation}
Let $\lambda = \frac12\min\{\lambda_1,\lambda_2\}>0$ and $\mathcal G_t = e^{\lambda t} \mathcal H_t$. Then $\mathcal G_0 = 0$, $\mathcal G_t\Ge 0$ and $\mathcal G_t$ satisfies
\begin{align}
\frac{\d}{\d t}\mathcal G_t & = 
\lambda e^{\lambda t} \mathcal F_t + e^{\lambda t}
\frac{\d}{\d t}\mathcal H_t \notag \\
& \Le \lambda e^{\lambda t} \mathcal H_t + 
e^{\lambda t} \Big(-\lambda_2 \mathcal H_t + C e^{-\lambda_3t}
\sqrt{\mathcal H_t} (|Q_0|+|P_0|)\delta\Big) \notag \\
& \Le -\lambda e^{\lambda t} \mathcal H_t + 
C \sqrt{\mathcal H_t} (|Q_0|+|P_0|) \delta \notag \\
& \Le - \lambda \mathcal G_t + C \sqrt{\mathcal G_t}
(|Q_0|+|P_0|)\delta.
\label{Gt decay}
\end{align}
Since $\mathcal G_0 = 0$, the inequality \eqref{Gt decay} implies
\begin{equation*}
	\sup_{t\Ge0} \sqrt{\mathcal G_t} \Le \frac{C}{\lambda}
	(|Q_0|+|P_0|)\delta.
\end{equation*}
As a consequence,
\begin{equation*}
	|W_t'-W_t| \Le C e^{-\frac12\lambda t}(|Q_0|+|P_0|)
	\delta,
\end{equation*}
which is exactly the desired result \eqref{qp diff decay}. \hfill $\blacksquare$\\[6pt]
\emph{Proof of Lemma~\ref{lemma: u2 estimate}. Estimate of the second derivative $\nabla^2 u(x,v,t)$.}\\[6pt]
According to \eqref{grad u expression}, $\nabla_x u(x,v,t)$ and $\nabla_v u(x,v,t)$ can be unified in the form
\begin{equation}
	g(x,v,t) = \mathbb E^{x,v}\Big[
	\nabla_x f(x_t,v_t)\cdot Q_t + 
	\nabla_v f(x_t,v_t)\cdot P_t
	\Big],
	\label{g xv}
\end{equation}
where the initial value $(Q_0,P_0)=(I_d,O_d)$ as in \eqref{Dx process} or $(Q_0,P_0)=(O_d,I_d)$ as in \eqref{Dv process}.

For another initial state $(x',v')$, let $(x_t',v_t')$ be the strong solution to the underdamped Langevin dynamics \eqref{underdamped} driven by the same Brownian motion $(B_t)_{t\Ge0}$. In other words, the stochastic processes $(x_t,v_t)_{t\Ge0}$ and $(x_t',v_t')_{t\Ge0}$ are in the synchronous coupling. Also, the tangent process $(Q_t',P_t')_{t\Ge0}$ is defined by
\begin{equation*}
(Q_t',P_t')_{t\Ge0} = \mathrm{TangentProcess}\big[
(x',v'),
(B_t)_{t\Ge0},(Q_0,P_0)
\big],
\end{equation*}
then the function $g(x',v',t)$ can be represented as
\begin{equation}
	g(x',v',t) = \mathbb E^{x',v'}
	\Big[
	\nabla_x f(x_t',v_t') \cdot Q_t' + \nabla_v 
	f(x_t',v_t')\cdot P_t'
	\Big].
	\label{g xv prime}
\end{equation}
Using \eqref{g xv} and \eqref{g xv prime}, the difference $g(x',v',t) - g(x,v,t)$ can be bounded by
\begin{align}
|g(x',v',t) - g(x,v,t)| \Le \,
& \mathbb E\big|
\nabla_x f(x_t',v_t')\cdot Q_t' - 
\nabla_x f(x_t,v_t)\cdot Q_t
\big| \notag \\ & \mathbb E\big|
\nabla_v f(x_t',v_t')\cdot P_t' - 
\nabla_v f(x_t,v_t)\cdot P_t
\big|.
\label{g difference}
\end{align}
Using the boundedness of $\nabla f(x,v)$ and $\nabla^2 f(x_t,v_t)$ in \textbf{A3}, we have
\begin{align}
& ~~~~ \mathbb E\big|
\nabla_x f(x_t',v_t')\cdot Q_t' - 
\nabla_x f(x_t,v_t)\cdot Q_t
\big| \notag \\
& \Le \mathbb E\big|\nabla_xf(x_t,v_t)\big| |Q_t' - Q_t| + \mathbb E\big|
\nabla_x f(x_t',v_t') - \nabla_x f(x_t,v_t)
\big||Q_t'| \notag \\
& \Le C\,\mathbb E\Big(|Q_t'-Q_t| + 
|x_t' - x_t| + |v_t' - v_t|
 \Big).
 \label{diff 1}
\end{align}
By Lemma~\ref{lemma: qp diff decay}, $|Q_t' - Q_t|$ has exponential decay in time, and for some constant $\lambda_1>0$,
\begin{equation}
|Q_t' - Q_t| \Le C e^{-\lambda_1 t} (|x' - x|+|v'-v|).
\label{diff 2}
\end{equation}
From the proof of Theorem~1 of \cite{coupling_1}, $|x_t'-x_t| + |v_t'-v_t|$ has exponential decay in time,
and for some constant $\lambda_2>0$,
\begin{equation}
|x_t' - x_t| + |v_t' - v_t| \Le C e^{-\lambda_2 t}
(|x'-x|+|v'-v|).
\label{diff 3}
\end{equation}
Let $\lambda = \min\{\lambda_1,\lambda_2\}>0$.
From the estimates \eqref{diff 1}\eqref{diff 2}\eqref{diff 3}, we obtain
\begin{equation*}
\mathbb E\big|
\nabla_x f(x_t',v_t')\cdot Q_t' - 
\nabla_x f(x_t,v_t)\cdot Q_t
\big| \Le C e^{-\lambda t}(|x'-x|+|v'-v|).
\end{equation*}
For the reason, it is easy to derive
\begin{equation*}
\mathbb E\big|
\nabla_v f(x_t',v_t')\cdot Q_t' - 
\nabla_v f(x_t,v_t)\cdot Q_t
\big| \Le C e^{-\lambda t}(|x'-x|+|v'-v|).
\end{equation*}
Therefore, \eqref{g difference} implies
\begin{equation*}
	|g(x',v',t) - g(x,v,t)| \Le Ce^{-\lambda t}(|x'-x|+|v'-v|).
\end{equation*}
As a consequence, the function $g(x,v,t)$ is global Lipschitz continuous in $(x,v)$ and
\begin{equation*}
|\nabla g(x,v,t)| \Le C e^{-\lambda t}.
\end{equation*}
Since $g(x,v,t)$ can be either $\nabla_x u(x,v,t)$ or $\nabla_v u(x,v,t)$, we conclude that 
\begin{equation*}
|\nabla^2 u(x,v,t)| \Le C e^{-\lambda t},
\end{equation*}
which is exactly the desired result. \hfill $\blacksquare$\\[6pt]
\emph{Proof of Lemma~\ref{lemma: SG strong order}. Unbiased property of stochastic gradient numerical integrators.}\\[6pt]
For simplicity, we focus on proving the unbiased condition \eqref{con 1} of \eqref{SG-UBU}. Using the explicit expressions in \eqref{UBU explicit}, we can represent $\tilde Z_1 = (\tilde X_1,\tilde V_1)$ and $\tilde Z_{0,1} = (\tilde X_{0,1},\tilde V_{0,1})$ as
\begin{align*}
\tilde Y_0 & = \tilde X_0 + \mathcal F(h/2) \tilde V_0 + 
\sqrt{2\gamma} \int_0^{h/2} \mathcal F(h/2-s)\d B_s, \\
\tilde X_1 & = X_0 + \mathcal F(h) V_0 - h\mathcal F(h/2)
b(\tilde Y_0,\omega) + \sqrt{2\gamma} \int_0^h 
\mathcal F(h-s)\d B_s, \\
\tilde V_1 & = \mathcal E(h)\tilde V_0 - h\mathcal E(h/2) 
b(\tilde Y_0,\omega) + \sqrt{2\gamma} \int_0^h 
\mathcal E(h-s)\d B_s, \\
\tilde X_{0,1} & = X_0 + \mathcal F(h) V_0 - h\mathcal F(h/2)
\nabla U(\tilde Y_0) + \sqrt{2\gamma} \int_0^h 
\mathcal F(h-s)\d B_s, \\
\tilde V_{0,1} & = \mathcal E(h)\tilde V_0 - h\mathcal E(h/2) 
\nabla U(\tilde Y_0) + \sqrt{2\gamma} \int_0^h 
\mathcal E(h-s)\d B_s.
\end{align*}
Then we have the explicit expressions of $\tilde X_1 - \tilde X_{0,1}$ and $\tilde V_1 - \tilde V_{0,1}$:
\begin{align*}
\tilde X_1 - \tilde X_{0,1} = h\mathcal F(h/2) \big(\nabla U(\tilde Y_0) - b(\tilde Y_0,\omega)\big),~~
\tilde V_1 - \tilde V_{0,1} = h\mathcal E(h/2)
\big(\nabla U(\tilde Y_0) - b(\tilde Y_0,\omega)\big).
\end{align*}
Since both $\tilde Y_0$ and $\tilde Z_{0,1} = (\tilde X_{0,1},\tilde V_{0,1})$ do not depend on $\omega$, we obtain
$$
\mathbb E\big[\tilde Z_1 - \tilde Z_{0,1}\big|\tilde Z_{0,1}\big]  = 0.
$$
Using the linear growth of $\nabla U(x)$ and $b(x,\omega)$ in \textbf{A1} and \textbf{B1}, we get
\begin{equation*}
	\mathbb E|\tilde Z_1 - \tilde Z_{0,1}|^4 \Le 
	Ch^4 \mathbb E(|\tilde Y_0|+1)^4 \Le Ch^4 
	\mathbb E(|\tilde Z_0|+1)^4,
\end{equation*}
which completes the proof. \hfill $\blacksquare$\\[6pt]
\emph{Proof of Lemma~\ref{lemma: moment SG}. Moment estimate for the stochastic gradient numerical integrator.}\\[6pt]
Note that $\tilde Z_0\{t\}$ is the solution to the underdamped Langevin dynamics
\begin{equation*}
\left\{
\begin{aligned}
\dot {\tilde X}_0\{t\} & = \tilde V_0\{t\}, \\
\dot {\tilde V}_0\{t\} & = - b(\tilde X_0\{t\},\omega) - 
\gamma \tilde V_0\{t\} + \sqrt{2\gamma} \dot B_t,
\end{aligned}
\right.
\end{equation*}
and $\mathcal H(x,v)$ defined in \eqref{Lyapunov} is a Lyapunov function, we have 
\begin{equation*}
	\mathbb E\big[\mathcal H_r(\tilde Z_0\{h\})\big]
	\Le e^{-ah} \mathcal H_r(\tilde Z_0) + \frac ba(1-e^{-ah}).
\end{equation*}
Since $\tilde Z_1$ is the time discretization approximation to $\tilde Z_0(h)$, utilizing the conditions \eqref{con 4} and \eqref{con 5}, similar to \eqref{Lyapunov derivation} we derive
\begin{align*}
	\mathbb E\big[\mathcal H_r(\tilde Z_1)\big] & \Le 
	\mathbb E\big[\mathcal H_r(\tilde Z_0\{h\})\big] + 
	\mathbb E\big|\mathcal H_r(\tilde Z_1) - \mathcal H_r(\tilde Z_0\{h\})\big|\\
	& \Le (e^{-ah} + Ch^{\frac32}) \mathcal H_r(\tilde Z_0) + 
	\frac{b}{a}(1-e^{-ah}),
\end{align*}
yielding the inequality
\begin{equation*}
	\mathbb E\big[\mathcal H_r(\tilde Z_1)\big]
	\Le e^{-a_1h} \mathcal H_r(\tilde Z_0) + 
	b_1 h
\end{equation*}
for some constants $a_1,b_1>0$. As a consequence, for the numerical solution $\tilde Z_n = (\tilde X_n,\tilde V_n)$,
\begin{equation*}
	\sup_{n\Ge0} \mathbb E\big[\mathcal H_r(\tilde Z_n)\big] \Le \max\bigg\{
	\frac{2b_1}{a_1} ,\mathcal H_r(\tilde Z_0)
	\bigg\} \Le C \mathbb E\big(|\tilde Z_0|+1\big)^{2r},
\end{equation*}
which gives rise to the uniform-in-time moment of the  numerical solution.
\hfill$\blacksquare$\\[6pt]
\emph{Proof of Lemma~\ref{lemma: diffusion SG}. Diffusion estimate for stochastic gradient numerical integrators.}\\[6pt]
In the proof of Lemma~\ref{lemma: diffusion}, the local error $\tilde Z_{0,1} - \tilde Z_0$ has been written as $\bar M_h + \bar N_h$, where the quantities $\bar M_h,\bar N_h$ satisfy
\begin{equation*}
\mathbb E\big[|\bar M_h|^2\big] \Le Ch (|\tilde Z_0|+1)^{2q},~~~~
\mathbb E\big[|\bar N_h|^2\big] \Le Ch^2 (|\tilde Z_0|+1)^{2q}.
\end{equation*}
Now let $\tilde M_h = \bar M_h + \tilde Z_1 - \tilde Z_{0,1}$, then we find that the local error $\tilde Z_1 - \tilde Z_0$ can be written as
\begin{equation*}
	\tilde Z_1 - \tilde Z_0 = \tilde M_h + \bar N_h,
\end{equation*}
where $\tilde M_h$ is a mean-zero variable and
\begin{equation*}
	\mathbb E\big[|\tilde M_h|^2\big] \Le 
	2 \mathbb E\big[|\bar M_h|^2\big] + 
	2 \mathbb E\big[|\tilde Z_1 - \tilde Z_{0,1}|^2\big]
	\Le Ch(|\tilde Z_0|+1)^{2q}.
\end{equation*}
The remaining procedure is the same with the proof of Lemma~\ref{lemma: diffusion}. \hfill $\blacksquare$
\end{document}